\documentclass[11pt]{article} 
\begin{document} 
\newtheorem{prop}{Proposition}[section]
\newtheorem{Def}{Definition}[section] \newtheorem{theorem}{Theorem}[section]
\newtheorem{lemma}{Lemma}[section] \newtheorem{Cor}{Corollary}[section]

\title{\bf The Cauchy problem for a Schr\"odinger - Korteweg - de Vries system 
with rough data}
\author{{\bf Hartmut Pecher}\\
Fachbereich Mathematik und Naturwissenschaften\\
Bergische Universit\"at Wuppertal\\
Gau{\ss}str.  20\\
D-42097 Wuppertal\\
Germany\\
e-mail Hartmut.Pecher@math.uni-wuppertal.de}
\date{}
\maketitle

\begin{abstract}
The Cauchy problem for a coupled system of the Schr\"odinger and the KdV 
equation is shown to be globally 
well-posed for data with infinite energy. The proof uses refined bilinear 
Strichartz type estimates and the 
I-method introduced by Colliander, Keel, Staffilani, Takaoka, and Tao.
\end{abstract}

\renewcommand{\thefootnote}{\fnsymbol{footnote}}
\footnotetext{\hspace{-1.8em}{\it 2000 Mathematics Subject Classification:} 
35Q53, 35Q55 }
\normalsize 
\setcounter{section}{-1}
\section{Introduction}
Consider the Cauchy problem 
\begin{eqnarray}
\label{1.1}
 iu_t + u_{xx} & = & \alpha vu + \beta |u|^2 u
\\ \label{1.2}
v_t + v v_x + v_{xxx} & = & \gamma (|u|^2)_x
\\
\label{1.3}
 u(x,0)  =  u_0(x)  & , & 
v(x,0)  =  v_0(x) 
\end{eqnarray}
where $u$ is a complex-valued und $v$ a
real-valued function defined for $(x,t) \in {\bf R}\times{\bf R}^+$ and 
$\alpha,\beta,\gamma \in {\bf R}$.

In the theory of capillary-gravity waves the interaction of a short and a long 
wave was modelled by such a 
coupled system of a Schr\"odinger and a KdV type equation (cf. Kawahara et al. 
\cite{KSK}, Funakoshi and 
Oikawa \cite{FO}). It also appears in plasma physics (cf. Nishikawa et al. 
\cite{NHMI}) modelling the 
interaction of the Langmuir and ion-acoustic waves.

The system (\ref{1.1}),(\ref{1.2}),(\ref{1.3}) was considered by M. Tsutsumi 
\cite{Ts} who showed that for 
$(u_0,v_0) \in H^{m+\frac{1}{2}}({\bf R}) \times H^m({\bf R})$ , $m=1,2,3,...$ , 
the problem is locally 
well-posed and for $\alpha \gamma > 0 $ also globally well-posed by using the 
conservation laws. Later, 
Bekiranov, Ogawa and Ponce \cite{BOP}, using the Fourier restriction norm 
method, lowered down the 
regularity assumptions on the data and proved local well-posedness for 
$(u_0,v_0) \in H^s({\bf R}) \times 
H^{s-\frac{1}{2}}({\bf R})$ for any $s \ge 0$.

Because our main aim is to consider the global problem which requires to use the 
conservation laws (cf. 
(\ref{3.0}),(\ref{3.0a}),(\ref{3.0b}) below) leading to an a-priori bound of the 
$H^1-$norms of $u$ and 
$v$, if $\alpha \gamma > 0$, we assume $u_0$ and $v_0$ to belong to the same 
Sobolev space $H^s({\bf R})$. 
Then we are able to show local well-posedness for any $s>0$ (cf. Theorem 
\ref{Theorem 1}) which implies 
especially 
global well-posedness in energy space $H^1({\bf R}) \times H^1{\bf R})$, if 
$\alpha \gamma > 0$. A global 
well-posedness result in this space was proven before by Guo and Miao \cite{GM} 
already. Moreover we are 
able to show global-wellposedness for less regular data, namely $u_0,v_0 \in 
H^s({\bf R})$ with $ s > 3/5 $ 
(if $\beta = 0$) and $ s > 2/3 $ (if $\beta \neq 0$). We use the Fourier 
restriction norm method and 
especially the I-method introduced by Colliander, Keel, Staffilani, Takaoka and 
Tao 
(\cite{CKSTT1},\cite{CKSTT},\cite{CKSTT2},\cite{CKSTT3},\cite{CKSTT4},
\cite{CKSTT5},\cite{CKSTT6}).
It was successfully 
applied to the (2+1)- and (3+1)-dimensional Schr\"odinger equation, to the 
(1+1)-dimensional derivative 
Schr\"odinger equation, and to the KdV and modified KdV equation with sometimes 
optimal results. In all 
these cases a scaling invariance was used which in our situation does not hold. 
Similar results were also 
given for the Klein - Gordon - Schr\"odinger system by Tzirakis \cite{Tz} and 
for the (1+1)-dimensional 
Zakharov system by the author \cite{P}.

The paper is organized as follows. In section 1 we give various bilinear 
Strichartz type estimates for the 
nonlinearities in the solution spaces $X^{s,b}$ and $Y^{s,b}$ for the 
Schr\"odinger and KdV equation, 
respectively, which are defined to be the completion of ${\cal S}({\bf 
R}\times{\bf R})$ w. r. to
$$ \|u\|_{X^{s,b}} := \|\langle \xi \rangle ^s \langle \tau + |\xi|^2 \rangle ^b 
\widehat{u}(\xi,\tau)\|_{L^2_{\xi \tau}} = \| \langle \xi \rangle ^s \langle 
\tau \rangle ^b {\cal 
F}(e^{-it\partial_x^2} u(x,t))\|_{L^2_{\xi \tau}} $$
and
$$ \|v\|_{Y^{s,b}} := \|\langle \xi \rangle ^s \langle \tau + \xi^3 \rangle ^b 
\widehat{v}(\xi,\tau)\|_{L^2_{\xi \tau}} = \| \langle \xi \rangle ^s \langle 
\tau \rangle ^b {\cal 
F}(e^{-t\partial_x^3} v(x,t))\|_{L^2_{\xi \tau}}. $$
For a given time interval $I$ we define $ \|u\|_{X^{s,b}(I)} := 
\inf_{\tilde{u}|I = u} 
\|\tilde{u}\|_{X^{s,b}}$ and similarly $ \|v\|_{Y^{s,b}(I)} $. The refined 
bilinear estimates are partly 
new or variants of known versions. In section 2 we formulate the local existence 
theorem and a variant of 
it for the modified system of differential equations after application of the 
operator $I$ to the original 
one (\ref{1.1}),(\ref{1.2}),(\ref{1.3}) giving a precise lower bound for the 
local existence time $T$ in 
terms of the norm of the data. This operator $I$, which gave the method its 
name, is defined as follows:
$I = I_N$ for given $s<1$ and $N >> 1$ is defined by $\widehat{I_N f}(\xi) := 
m_N(\xi) \widehat{f}(\xi)$. 
Here $m_N(\xi)$ is a smooth, radially symmetric and nonincreasing function of 
$|\xi|$, defined by $m_N(\xi) 
= 1$ for $|\xi| \le N$ and $ m_N(\xi) = (\frac{N}{|\xi|})^{1-s}$ for $|\xi| \ge 
2N$. Dropping $N$ from the 
notation we have $ I: H^s \to H^1 $ is a smoothing operator with
$$ \|u\|_{X^{m,b}} \le c \|Iu\|_{X^{m+1-s,b}} \le c N^{1-s} \|u\|_{X^{m,b}} $$
and similarly for $Y^{s,b}$. In section 3 we consider the conserved quantities 
$L(u,v)$ and $E(u,v)$ 
(cf. (\ref{3.0a}),(\ref{3.0b})) and their modified versions $L(Iu,Iv)$ and 
$E(Iu,Iv)$ which are no longer 
conserved, but it is possible to control their growth in time, because here some 
sort of cancellation 
helps. As is typical for the I-method we then in sections 4 and 5 consider in 
detail this increment of 
$E(Iu,Iv)$ and $L(Iu,Iv)$ , respectively, which is shown to be small for small 
time intervals and large 
$N$. Important tools here are the refined bilinear Strichartz estimates of 
section 1, especially a new 
estimate for the product of a Schr\"odinger and a KdV part. These estimates also 
allow to control the 
growth 
of the corresponding norms of the solution during its time evolution. One 
iterates in section 6 the local 
existence theorem with time steps of equal length. To achieve this one has to 
make the process uniform 
which can be done if $s$ is close enough to $1$.\\ 
We collect some elementary facts about the spaces $X^{m,b}$ (and analogously 
$Y^{m,b}$).\\ 
The following interpolation property is well-known:\\
$ X^{(1-\Theta)m_0+\Theta m_1,(1-\Theta)b_0+\Theta b_1} = 
(X^{m_0,b_0},X^{m_1,b_1})_{[\Theta]} $ for $ \Theta \in [0,1] $. \\
If $u$ is a 
solution of $iu_t+\partial_x^2 u = 0$ with $u(0)=f$ and $\psi$ is a cutoff 
function in 
$C^{\infty}_0({\bf R})$ with $supp \, \psi \subset (-2,2)$ , $\psi \equiv 1$ on 
$[-1,1]$ , $ \psi(t) = 
\psi(-t) $ , $ \psi(t)\ge 0 $ , $\psi_{\delta}(t):=\psi(\frac{t}{\delta}) \, , $ 
 $ 0<\delta \le 1$, we 
have 
for $b>0$: 
$$\|\psi_1 u\|_{X^{m,b}} \le c \|f\|_{H^m} \, . $$
If $v$ is a solution of the problem $iv_t +\partial_x^2 v = F $ , $ v(0)=0 $ , 
we have for $b'+1 
\ge b \ge 0 \ge b' > -1/2$
$$ \|\psi_{\delta} v \|_{X^{m,b}} \le c \delta^{1+b'-b} \|F\|_{X^{m,b'}} $$
(for a proof cf. \cite{GTV}, Lemma 2.1).\\
Finally, if $1/2 > b > b' \ge 0$ , $ m \in {\bf R} $ , we have the embedding
\begin{equation}
\label{0}
\|f\|_{X^{m,b'}[0,\delta]} \le c \delta^{b-b'} \|f\|_{X^{m,b}[0,\delta]} \, .
\end{equation}
For the convenience of the reader we repeat the proof of \cite{G}, Lemma 1.10. 
The claimed estimate is an 
immediate consequence of the following
\begin{lemma}
For $1/2 > b > b' \ge 0$ , $ 0 < \delta \le 1 $ , $ m \in {\bf R} $ the 
following estimate holds:
$$ \| \psi_{\delta} f \|_{X^{m,b'}} \le c \delta^{b-b'} \|f\|_{X^{m,b}} \, . $$
\end{lemma}
{\bf Proof:} The following Sobolev multiplication rule holds:
$$ \|fg\|_{H^{b'}_t} \le c \|f\|_{H^{\frac{1}{2}-(b-b')}_t} \|g\|_{H^b_t} \, . 
$$
This rule follows easily by the Leibniz rule for fractional derivatives, using 
$J^s := {\cal 
F}^{-1}\langle\tau\rangle^s {\cal F}$:
\begin{eqnarray*}
\|fg\|_{H^{b'}_t} & \le & c(\|(J^{b'}f)g\|_{L^2_t} + \|f(J^{b'}g)\|_{L^2_t}) \\
& \le & c(\|J^{b'}f\|_{L^p_t} \|g\|_{L^{p'}_t} + \|f\|_{L^{q'}_t} 
\|J^{b'}g\|_{L^q_t})
\end{eqnarray*}
with $\frac{1}{p}=b$ , $ \frac{1}{p'} = \frac{1}{2}-b $ , $ \frac{1}{q'}=b-b' $ 
, $ \frac{1}{q} = 
\frac{1}{2} - (b-b') $. Sobolev's embedding theorem gives the claimed result. 
Consequently we get
$$ \|\psi_{\delta}g\|_{H^{b'}_t} \le c 
\|\psi_{\delta}\|_{H^{\frac{1}{2}-(b-b')}_t} \|g\|_{H^b_t} \le c 
\delta^{b-b'} \|g\|_{H^b_t} \, , $$
and thus
\begin{eqnarray*}
\|\psi_{\delta}f\|_{X^{m,b'}} = 
\|e^{-it\partial_x^2}\psi_{\delta}f\|_{H^m_x\otimes 
H^{b'}_t} & \le & c \delta^{b-b'}\|e^{-it\partial_x^2}f\|_{H^m_x \otimes 
H^b_t}\\ 
& = & c\delta^{b-b'}\|f\|_{X^{m,b}_{\varphi}} \, . 
\end{eqnarray*}
Fundamental are the following linear Strichartz type estimates for the 
Schr\"odinger equation (cf. e.g. 
\cite{GTV}, Lemma 2.4):
$$ \|e^{it\partial^2_x} \psi\|_{L^q_t(I,L_x^r({\bf R}))} \le c 
\|\psi\|_{L^2_x({\bf R})} $$
and
$$ \|u\|_{L^q_t(I,L_x^r({\bf R}))} \le c \|u\|_{X^{0,\frac{1}{2}+}(I)} \, , $$
if $ 0 \le \frac{2}{q} = \frac{1}{2} - \frac{1}{r} $ , especially
$$ \|u\|_{L^6_{xt}} \le c \|u\|_{X^{0,\frac{1}{2}+}} \, , $$
which by interpolation with the trivial case $\|u\|_{L^2_{xt}} = 
\|u\|_{X^{0,0}}$ gives:
$$ \|u\|_{L^p_{xt}} \le c \|u\|_{X^{0,\frac{3}{2}(\frac{1}{2}-\frac{1}{p})+}} \, 
, $$
if $ 2 < p \le 6 $.\\
For the KdV (Airy) equation we have (cf. e.g. \cite{KPV1}, Theorem 2.4): \\ $ 
\|e^{-t\partial_x^3} 
\psi\|_{L^8_{xt}} \le c 
\|\psi\|_{L^2_x} $ , and thus $ \|v\|_{L^8_{xt}} \le c 
\|v\|_{Y^{0,\frac{1}{2}+}}$.\\
We use the notation $\langle \lambda \rangle := (1+\lambda^2)^{1/2} $. Let $a 
\pm$ denote a number 
slightly larger (resp., smaller) than $a$.
\section{Bilinear Strichartz type estimates}
\begin{lemma}
\label{Lemma 1}
If $ s \ge 0 $ , $ b > \frac{1}{2} $ , $ b' \ge max(\frac{1}{4}-\frac{s}{3},0) $ 
, the following estimate 
holds:
\begin{equation}
\label{1.1'}
\| \partial_x (v_1 v_2) \|_{Y^{s,-b'}} \le c \|v_1\|_{Y^{s,b}} \|v_2\|_{Y^{s,b}} 
\, .
\end{equation} 
\end{lemma}
{\bf Proof:} With $\xi_i \in supp \, \widehat{v_i} $ $(i=1,2)$ we first consider 
the case $|\xi_1| >> 
|\xi_2|$ (or similarly $|\xi_2| >> |\xi_1|$). In this case $|\xi|:=|\xi_1 
+\xi_2| \sim |\xi_1| \sim |\xi_1 
+ \xi_2|^{\frac{1}{2}} |\xi_1 - \xi_2|^{\frac{1}{2}}$, and we use (\ref{g}) to 
conclude
$$ \| \partial_x (v_1 v_2)\|_{L^2_{xt}} \le c \|v_1\|_{Y^{0,b}} 
\|v_2\|_{Y^{0,b}} \, , $$
which immediately implies
$$ \| \partial_x (v_1 v_2)\|_{Y^{s,0}} \le c \|v_1\|_{Y^{s,b}} \|v_2\|_{Y^{s,b}} 
\, . $$
Next we consider the case $|\xi_1| \sim |\xi_2|$. This implies $ |\xi| \le 
|\xi_1|,|\xi_2|$. We have to 
show
$$ \left\|\frac{\xi \langle \xi \rangle ^s}{\langle \tau - \xi ^3 \rangle ^{b'}} 
\int \int 
\frac{f_1(\xi_1,\tau_1)}{\langle \tau_1 - \xi_1^3 \rangle ^b \langle \xi_1 
\rangle ^s} \cdot 
\frac{f_2(\xi_2,\tau_2)}{\langle \tau_2 - \xi_2^3 \rangle ^b \langle \xi_2 
\rangle ^s} d\xi_1 d\tau_1 
\right\|_{L^2_{\xi \tau}} \le c \|f_1\|_{L^2_{\xi \tau}} \|f_2\|_{L^2_{\xi 
\tau}} \, , $$
where $ \xi = \xi_1 + \xi_2$ , $ \tau = \tau_1 + \tau_2$. Using the Schwarz 
method this is implied by
$$ \sup_{\xi,\tau} \frac{ |\xi| \langle \xi \rangle ^s}{\langle \tau - \xi ^3 
\rangle ^{b'}} \left( \int 
\int \frac{d\tau_1 d\xi_1}{ \langle \tau_1 - \xi_1 ^3 \rangle ^{2b} \langle 
\tau_2 - \xi_2 ^3 \rangle ^{2b} 
\langle \xi_1 \rangle ^{2s} \langle \xi_2 \rangle ^{2s} } \right)^{\frac{1}{2}} 
\le c \, . $$
Using (the proof of) \cite{KPV2}, Lemma 2.4, the l.h.s. is bounded by
\begin{eqnarray*}
\lefteqn{c \sup_{\xi,\tau} \frac{|\xi|}{\langle \xi\rangle ^s \langle \tau - \xi 
^3 \rangle ^{b'}} 
\left(\int \int \frac{d\tau_1 \, d\xi_1}{\langle \tau_1 - \xi_1^3 \rangle ^{2b} 
\langle \tau_2 - \xi_2^3 
\rangle ^{2b}} \right)^{\frac{1}{2}}} \\
& \le & c \sup_{\xi,\tau} \frac{|\xi|}{\langle \xi\rangle ^s \langle \tau - \xi 
^3 \rangle 
^{b'}|\xi|^{\frac{1}{4}} \langle 4\tau - \xi^3 \rangle ^{\frac{1}{4}}} \, .
\end{eqnarray*}
This can easily seen to be finite, if $ \frac{3}{4} -s-3b' \le 0 $ and $ b' \ge 
0 $ , which is assumed.
\begin{lemma}
\label{Lemma 2} 1. If $s \ge 0$ , $ b > 1/2 $ , $ b' > 
\max(\frac{1}{6},\frac{1}{2}-s) $ , we have
\begin{equation}
\label{1.2'}
\| \partial_x (u_1 \overline{u_2}) \|_{Y^{s,-b'}} \le c \|u_1\|_{X^{s,b}} 
\|u_2\|_{X^{s,b}} \, .
\end{equation}
2. If  $ b,b' > 1/2 $ , we have
\begin{equation}
\label{1.3'}
\| \partial_x (u_1 \overline{u_2}) \|_{Y^{0,-b'}} \le c \|u_1\|_{X^{0,b}} 
\|u_2\|_{X^{0,b}} \, .
\end{equation}
\end{lemma}
{\bf Proof:} First we consider 2., which is proved along the lines of 
\cite{BOP}, Lemma 3.2. With the 
restriction $ \xi = \xi_1 + \xi_2 $ , $ \tau = \tau_1 + \tau_2 $ , $ \sigma = 
\tau - \xi^3 $ , $ \sigma_1 = 
\tau_1 - \xi_1^2 $ , $ \sigma_2 = \tau_2 + \xi_2^2 $ , we have to show
$$ \left| \int \frac{|\xi| g(\xi,\tau)}{\langle \sigma \rangle ^{b'}} \cdot 
\frac{f_1(\xi_1,\tau_1)}{\langle \sigma_1 \rangle ^b}\cdot 
\frac{f_2(\xi_2,\tau_2)}{\langle \sigma_2 
\rangle ^b} \, d\tau_1 \, d\xi_1 \, d\tau \, d\xi \right|  
 \le  c \|g\|_{L^2_{xt}}\|f_1\|_{L^2_{xt}}\|f_2\|_{L^2_{xt}} \, . $$
Using Schwarz' method this is implied, if one of the following estimates holds 
in each of the regions into 
which we are going to split our domain:
\begin{eqnarray*}
\sup_{\tau,\xi} \frac{|\xi|}{\langle \sigma \rangle ^{b'}} \left( \int \int 
\frac{d\xi_1 \, 
d\tau_1}{\langle \sigma_1 \rangle ^{2b} \langle \sigma_2 \rangle ^{2b}} 
\right)^{\frac{1}{2}} & \le & c \, 
, \\
\sup_{\tau_1,\xi_1} \frac{1}{\langle \sigma_1 \rangle ^{b}} \left( \int \int 
\frac{|\xi|^2 \, d\xi \, 
d\tau}{\langle \sigma \rangle ^{2b'} \langle \sigma_2 \rangle ^{2b}} 
\right)^{\frac{1}{2}} & \le & c \, , 
\\ 
\sup_{\tau_2,\xi_2} \frac{1}{\langle \sigma_2 \rangle ^{b}} \left( \int \int 
\frac{|\xi|^2 \, d\xi_1 \, 
d\tau_1}{\langle \sigma \rangle ^{2b'} \langle \sigma_1 \rangle ^{2b}} 
\right)^{\frac{1}{2}} & \le & c \, . 
\end{eqnarray*} 
Here the integrations w.r. to $\tau$ and $\tau_1$ can be done using \cite{BOP}, 
Lemma 2.5, and lead to the 
following sufficient conditions, using $b>1/2$ and w.l.o.g. $b' \le b$ :
\begin{eqnarray}
\label{i}
\sup_{\tau,\xi} \frac{|\xi|}{\langle \sigma \rangle ^{b'}} \left( \int 
\frac{d\xi_1}{\langle \tau 
+(\xi-\xi_1)^2-\xi_1^2 \rangle ^{2b}} \right)^{\frac{1}{2}} & \le & c \, , \\
\label{ii}
\sup_{\tau_1,\xi_1} \frac{1}{\langle \sigma_1 \rangle ^{b}} \left( \int 
\frac{|\xi|^2 \,d\xi}{\langle \xi^3 
- \tau_1 +(\xi-\xi_1)^2 \rangle ^{2b'}} \right)^{\frac{1}{2}} & \le & c \, ,  \\
\label{iii}
\sup_{\tau_2,\xi_2} \frac{1}{\langle \sigma_2 \rangle ^{b}} \left( \int 
\frac{|\xi|^2 \, d\xi_1}{\langle 
\xi_1^2 + \tau_2 -(\xi_1 + \xi_2)^3 \rangle ^{2b'}} \right)^{\frac{1}{2}} & \le 
& c \, . 
\end{eqnarray}
Now we split our domain $(\xi,\xi_1,\tau,\tau_1) \in {\bf R}^4$ as follows: 
${\bf R}^4 = A \cup B \cup C $ 
, where $A=\{ |\xi| \le 10\} $ , $ B = \{ |\xi| > 10 \, , \, 
|3\xi^2-2(\xi-\xi_1)| \ge \xi^2 \} $ , $ C = 
\{ |\xi| > 10 \, , \, |\xi^2+\xi-2\xi_1| \ge \frac{1}{2}\xi^2 \} $ .\\
 We remark here (for later reference) that the condition in $C$ is especially 
fulfilled, if $ |\xi_2| >> 
|\xi_1| $ (or $ |\xi_1| >> |\xi_2| $) and $|\xi| > 10$.\\
That we in fact have ${\bf R}^4 = A \cup B \cup C $ can be seen as follows: if 
both conditions defining 
$B$ and $C$ are violated, i.e. $ |3\xi^2-2(\xi-\xi_1)| < \xi^2 $ and $ 
|\xi^2+\xi-2\xi_1| < 
\frac{1}{2}\xi^2 ,$ then we have
$$ |2\xi_1| = |(2\xi_1-2\xi+3\xi^2)+2\xi-3\xi^2| > 3\xi^2 -2|\xi|-\xi^2 = 2 
\xi^2 - 2|\xi| $$
and
$$ |2\xi_1| = |(2\xi_1-\xi -\xi^2)+\xi +\xi^2| < \frac{1}{2} \xi^2 +|\xi| +\xi^2 
= \frac{3}{2} \xi^2 +|\xi| 
\, . $$
Thus
$$ 2\xi^2-2|\xi| < \frac{3}{2} \xi^2 + |\xi| \Leftrightarrow \frac{1}{2}\xi^2 < 
3|\xi| \Leftrightarrow 
|\xi| < 6 \, , $$
thus we are in $A$ in this case.\\
For the region $A$ we trivially have (\ref{ii}).\\
For the region $B$ we again prove (\ref{ii}). The change of variables $\eta := 
\xi^3-\tau_1+(\xi-\xi_1)^2$ 
gives $ d\xi = \frac{d\eta}{3\xi^2+2(\xi-\xi_1)} $. Using the definition of $B$ 
we thus estimate the l.h.s. 
of (\ref{ii}) by
$$ \sup_{\tau_1,\xi_1} \frac{1}{\langle \tau_1 - \xi_1^2 \rangle ^b} \left( \int 
\frac{d\eta}{\langle \eta 
\rangle ^{2b'}} \right)^{\frac{1}{2}} \le c \, , $$
provided $ b \ge 0$ , $ b' > 1/2 $ .\\
For the region $C$ we use the algebraic inequality
\begin{eqnarray*}
|\sigma| + |\sigma_1| + |\sigma_2| & \ge & 
|(\tau-\xi^3)-(\tau_1-\xi_1^2)-(\tau_2+\xi_2^2)| \\
& = & |-\xi^3+\xi_1^2-\xi_2^2| = |-\xi^3+\xi_1^2-(\xi-\xi_1)^2| \\
& = & |-\xi^3-\xi^2+2\xi \xi_1| = |\xi||\xi^2+\xi-2\xi_1| \ge \frac{1}{2} 
|\xi|^3 \, ,
\end{eqnarray*}
which leads to 3 cases, depending on which one of the 3 terms on the l.h.s. is 
dominant.\\
If $|\sigma|$ is dominant, we prove (\ref{i}). Its l.h.s. is estimated for 
$b>1/2$ and $b' > 1/6$ by
$$ c \sup_{\tau,\xi} |\xi|^{1-3b'} \left( \int \frac{d\xi_1}{\langle \tau + 
\xi^2 -2\xi \xi_1 \rangle 
^{2b}} \right)^{\frac{1}{2}} \le c \sup_{\xi} |\xi|^{\frac{1}{2}-3b'} \le c \, . 
$$
If $|\sigma_1|$ is dominant, we estimate the l.h.s. of (\ref{ii}) for $b> 1/2$ 
and $b' \ge 0$ by
$$ c \sup_{\tau_1,\xi_1} \left( \int \frac{d\xi}{|\xi|^{6b-2} \langle 
\xi^3-\tau_1+(\xi-\xi_1)^2 \rangle 
^{2b'}}\right)^{\frac{1}{2}} \le c \, . $$
If $|\sigma_2|$ is dominant, we estimate the l.h.s. of (\ref{iii}) by 
\cite{BOP}, Lemma 2.5 for $b> 1/2$ 
and 
$b' > 1/6$ :
$$ c \sup_{\tau_2,\xi_2} \left( \int \frac{d\xi_1}{\langle \xi \rangle ^{6b-2} 
\langle \xi_1^2 + \tau_2 
-(\xi_1 + \xi_2)^3 \rangle ^{2b'}} \right)^{\frac{1}{2}} \le c \, .  $$
This proves claim 2.\\
Next we prove claim 1. In the regions $|\xi_1| >> |\xi_2| $ or $|\xi_2| >> 
|\xi_1| $ we can immediately use 
our considerations above for the regions $A$ and $C$ and get for $b' > 1/6$ , $b 
> 1/2$ :
$$ \|\partial_x(u_1 \overline{u_2})\|_{Y^{0,-b'}} \le c \|u_1\|_{X^{0,b}} 
\|u_2\|_{X^{0,b}} \, ,  $$
and consequently for $s \ge 0$ :
$$ \|\partial_x(u_1 \overline{u_2})\|_{Y^{s,-b'}} \le c \|u_1\|_{X^{s,b}} 
\|u_2\|_{X^{s,b}} \, . $$
In the region $|\xi_1| \sim |\xi_2|$ we interpolate our estimate 2. with the 
following estimate of 
\cite{BOP}, Lemma 3.2:
\begin{equation}
\label{*}
\|\partial_x(u_1 \overline{u_2})\|_{Y^{-\frac{1}{2},0}} \le c \|u_1\|_{X^{0,b}} 
\|u_2\|_{X^{0,b}} \, ,
\end{equation}
and get for $ 0 \le l \le 1/2 $ :
$$ \|\partial_x(u_1 \overline{u_2})\|_{Y^{-l,l-\frac{1}{2}-}} \le c 
\|u_1\|_{X^{0,b}} \|u_2\|_{X^{0,b}} \, 
. $$
This implies, using $|\xi| \le |\xi_1|,|\xi_2|$ and $b'>\frac{1}{2}-s$, for $0 
\le s \le \frac{1}{2}$ :
\begin{eqnarray*}
\hspace{-1cm} \|\partial_x(u_1 \overline{u_2})\|_{Y^{s,-b'}} & \le & c \| J^{2s} 
\partial_x(u_1 
\overline{u_2})\|_{Y^{-s,-b'}} \\
\le c \|\partial_x(J^s u_1 J^s \overline{u_2})\|_{Y^{-s,-b'}} & \le & c 
\|\partial_x(J^s u_1 J^s 
\overline{u_2})\|_{Y^{-s,s-\frac{1}{2}-}} \\
\le \; c \|J^s u_1\|_{X^{0,b}} \|J^s u_2\|_{X^{0,b}} & \le & c 
\|u_1\|_{X^{s,b}}\|u_2\|_{X^{s,b}} \, .
\end{eqnarray*}
If $s > 1/2$ we estimate similarly for $b'=0$ using (\ref{*}).\\
This completes the proof of Lemma \ref{Lemma 2}.
\begin{lemma}
\label{Lemma 3}
If $\widehat{u}(\xi_1,t)$ and $\widehat{v}(\xi_2,t)$ are supported in the region 
$|\xi_2|^2 >> |\xi_1|$ , 
we have the following inequality for $b > 1/2$ :
\begin{equation}
\label{1.4}
\|u D_x v\|_{L^2_{xt}} \le c \|u\|_{X^{0,b}} \|v\|_{Y^{0,b}} \, .
\end{equation}
Here $u$ on the l.h.s. can be replaced by $\overline{u}$.
\end{lemma}
{\bf Proof:} a) If $|\xi_2|$ is bounded, this follows directly from Strichartz' 
estimates:
$$ \|u D_x v\|_{L^2_{xt}} \le c \|u\|_{L^4_{xt}} \|D_x v\|_{L^4_{xt}} \le c 
\|u\|_{X^{0,b}} \|D_x 
v\|_{Y^{0,b}} \le c \|u\|_{X^{0,b}} \|v\|_{Y^{0,b}} \, . $$
b) It suffices to prove for functions $\widehat{u_0}(\xi_1)$ and $ 
\widehat{v_0}(\xi_2)$ with support in $ 
\{ |\xi_2|^2 >> |\xi_1| \, , \, |\xi_2| >> 1 \} $ :
$$ \| e^{it\partial_x^2} u_0 e^{-t\partial_x^3} D_x v_0\|_{L^2_{xt}} \le c 
\|u_0\|_{L^2_x} \|v_0\|_{L^2_x}. 
$$
We have
\begin{eqnarray*}
\lefteqn{ \| e^{it\partial_x^2} u_0 e^{-t\partial_x^3} D_x v_0\|_{L^2_{xt}} } \\
& = & \int d\xi dt \int_* d\xi_2 d\eta_2 \, e^{-it(\xi_1^2 + \xi_2^3 - \eta_1^2 
- \eta_2^3)} |\xi_2 \eta_2| 
\widehat{u_0}(\xi_1) \overline{\widehat{u_0}(\eta_1)} \widehat{v_0}(\xi_2) 
\overline{\widehat{v_0}(\eta_2)} 
\\
& = & \int d\xi \int_* d\xi_2 d\eta_2 \, \delta(P(\eta_2)) |\xi_2 \eta_2| 
\widehat{u_0}(\xi_1) 
\overline{\widehat{u_0}(\eta_1)} \widehat{v_0}(\xi_2) 
\overline{\widehat{v_0}(\eta_2)} \, .
\end{eqnarray*}
Here * denotes integration over $\xi = \xi_1 + \xi_2 = \eta_1 + \eta_2 $ , $ 
|\xi_2|^2 >> |\xi_1| , $  $ 
|\eta_2|^2 >> |\eta_1|$ , $ |\xi_2|,|\eta_2| >> 1$ , and
\begin{eqnarray*}
P(\eta_2) & := & \eta_1^2 + \eta_2^3 - \xi_1^2 - \xi_2^3 = (\xi - \eta_2)^2 + 
\eta_2^3 - (\xi - \xi_2)^2 - 
\xi_2^3 \\
& = & (\eta_2 -\xi_2)(\eta_2^2 + \eta_2(1+\xi_2)-2\xi + \xi_2 + \xi_2^2).
\end{eqnarray*}
Now $P(\eta_2)$ has the root $\eta_2 = \xi_2$ and no further real zero, because
$$ \eta_2^2+\eta_2(1+\xi_2)-2\xi+\xi_2 + \xi_2^2 = (\eta_2 + 
\frac{1+\xi_2}{2})^2 + \frac{3}{4} \xi_2^2 - 
\frac{3}{2}\xi_2 -2\xi_1 - \frac{1}{4} > 0 $$
for $|\xi_2|^2 >> |\xi_1| $ , $ |\xi_2| >> 1 $ . Moreover $P'(\eta_2) = 
3\eta_2^2 -2(\xi - \eta_2)$ . Thus 
for $|\xi_2|^2 >> |\xi_1| $ we get $ |P'(\xi_2)| = |3\xi_2^2-2(\xi - \xi_2)| = 
|3\xi_2^2 - 2\xi_1| \ge 
\xi_2^2 $. Now we use the well-known identity $\delta(P(\eta_2)) = 
\frac{\delta(\eta_2 - 
\xi_2)}{|P'(\xi_2)|} $ and estimate as follows:
\begin{eqnarray*}
\lefteqn{\| e^{it\partial_x^2} u_0 e^{-t\partial_x^3} D_x v_0\|_{L^2_{xt}}^2} \\ 
& &\le  c \int d\xi \int_* 
d\xi_2 d\eta_2 \frac{\delta(\eta_2-\xi_2)}{\xi_2^2} |\xi_2 \eta_2| 
\widehat{u_0}(\xi_1) 
\overline{\widehat{u_0}(\eta_1)} \widehat{v_0}(\xi_2) 
\overline{\widehat{v_0}(\eta_2)} 
 =  c \|u_0\|_{L^2_x}^2 \|v_0\|_{L^2_x}^2.
\end{eqnarray*}
\\[0.5cm]
Further bilinear and trilinear estimates used in the sequel are the following:
\begin{equation}
\label{a}
\|uv\|_{X^{s,0}} \le c \|u\|_{X^{s,b}} \|v\|_{Y^{s,b}}
\end{equation}
for $s\ge 0 \, , \, b > 1/2 $ , which follows easily for $s=0$ by Strichartz' 
estimates, and thus for any 
$s \ge 0$. This holds true for $u$ replaced by $\overline{u}$ on the l.h.s.\\
Similarly one gets:
\begin{equation}
\label{b}
\| |u|^2 u\|_{X^{s,0}} \le c \|u\|_{X^{s,b}}^3
\end{equation}
for $s \ge 0$ , $ b > 1/2 $ .\\
Defining
$$ {\cal F} I_-^{\alpha} (u,v)(\xi,\tau) := \int_* d\xi_1 d\tau_1 \, |\xi_1 - 
\xi_2|^{\alpha} 
\widehat{u}(\xi_1,\tau_1) \widehat{v}(\xi_2,\tau_2) \, ,  $$
where * denotes integration over $\xi = \xi_1 + \xi_2$ , $ \tau = \tau_1 + 
\tau_2$ , we have by \cite{G1}, 
Corollary 3.2 for $b > 1/2$ :
\begin{equation}
\label{g}
\|D_x^{\frac{1}{2}} I_-^{\frac{1}{2}}(v_1,v_2)\|_{L^2_{xt}} \le c 
\|v_1\|_{Y^{0.b}} \|v_2\|_{Y^{0,b}}\, . 
\end{equation}
Moreover we have 
\begin{equation}
\label{h}
\|I_-^{\frac{1}{2}}(u_1,u_2)\|_{L^2_{xt}} \le c \|u_1\|_{X^{0,b}} 
\|u_2\|_{X^{0.b}}
\end{equation}
for $ b > 1/2 $ , cf. \cite{G}, Lemma 4.2 or \cite{CKSTT1}, (the proof of) Lemma 
7.1, and especially
\begin{equation}
\label{c}
\|(D_x^{\frac{1}{2}} u_1)u_2\|_{L^2_{xt}} \le c \|u_1\|_{X^{0,b}} 
\|u_2\|_{X^{0,b}} \, ,
\end{equation}
if $ |\xi_1| \ge \beta |\xi_2| $ for $\xi_j \in supp \, \widehat{u_j} $ 
$(j=1,2)$ , where $\beta > 1 $ .
The last inequality remains true for $u_1$ and/or $u_2$ replaced by 
$\overline{u_1}$ and/or 
$\overline{u_2}$ on the l.h.s. Also
\begin{equation}
\label{d}
\|D_x^{\frac{1}{2}} (u_1 \overline{u_2})\|_{L^2_{xt}} \le c \|u_1\|_{X^{0,b}} 
\|u_2\|_{X^{0,b}}
\end{equation}
for $ b > 1/2 $ (cf. \cite{BOP2}, Lemma 3.2), and
\begin{equation}
\label{e}
\|v_1 v_2\|_{L^2_{xt}} \le c \|v_1\|_{Y^{-\frac{1}{2},\frac{1}{2}-}} 
\|v_2\|_{Y^{\frac{1}{4},\frac{1}{2}+}}
\end{equation}
(cf. \cite{T}, Prop. 6.2). Finally
\begin{equation}
\label{f}
\|\partial_x(v_1 v_2)\|_{X^{0,-\frac{1}{2}+}} \le c 
\|v_1\|_{Y^{\gamma_1,\frac{1}{2}+}} 
\|v_2\|_{Y^{\gamma_2,\frac{1}{2}+}}
\end{equation}
if $\widehat{v_i}$ are supported outside $|\xi| \le 1$ and $ \gamma_1 + \gamma_2 
> - \frac{3}{4} $ ; $ 
\gamma_1 , \gamma_2 > - \frac{1}{2} $ (cf. \cite{CKSTT}, Lemma 1).
\section{Local well-posedness}
Consider the Cauchy problem (\ref{1.1}),(\ref{1.2}),(\ref{1.3}), where 
$\alpha,\beta,\gamma \in {\bf R}$. 
Using the Fourier restriction norm method it is not difficult to prove the 
following local well-posedness 
result using the estimates in chapter 1.
\begin{theorem}
\label{Theorem 1}
For any $(u_0,v_0) \in H^s({\bf R}) \times H^s({\bf R})$ and $ s > 0 $ there 
exists $ b > 1/2 $ and 
$\delta=\delta(\|u_0\|_{H^s},\|v_0\|_{H^s}) > 0 $ , such that 
(\ref{1.1}),(\ref{1.2}),(\ref{1.3}) has a 
unique solution satisfying $(u,v) \in X^{s,b}[0,\delta] \times Y^{s,b}[0,\delta] 
$ and $(u,v) \in 
C^0([0,\delta],H^s({\bf R}) \times H^s({\bf R}))$ . This solution depends 
continuously on the data 
$(u_0,v_0)$ .
\end{theorem}
{\bf Proof:} One constructs a fixed point of the mapping $S=(S_0,S_1)$ induced 
by the corresponding 
integral equations:
\begin{eqnarray*}
S_0(u,v) & := & e^{it\partial_x^2} u_0 - i \int_0^t 
e^{i(t-s)\partial_x^2}[\alpha u(s)v(s) + \beta |u(s)|^2 
u(s)] \, ds \, , \\
S_1(u,v) & := & e^{-t\partial_x^3} v_0 - \int_0^t e^{-(t-s)\partial_x^3} 
[v(s)\partial_x v(s) - \gamma 
\partial_x(|u(s)|^2)] \, ds \, .
\end{eqnarray*}
In order to estimate the nonlinear terms we use (\ref{a}),(\ref{b}), Lemma 
\ref{Lemma 1} and Lemma 
\ref{Lemma 2}, which allows to choose $b' < 1/2$, if $ s > 0$ . Standard 
arguments using some of the facts 
given in the introduction imply the claimed result. 

In order to apply the I-method we also need a variant of this local result. 
Applying the I-operator to the 
system (\ref{1.1}),(\ref{1.2}),(\ref{1.3}) we get
\begin{eqnarray}
\label{1.1''}
iIu_t + Iu_{xx} & = & \alpha I(uv) + \beta I(|u|^2u)\, ,  \\
\label{1.2''}
Iv_t + I(vv_x) + I(v_{xxx}) & = & \gamma I(|u|^2)_x \, ,\\
\label{1.3''}
Iu(0) = Iu_0 & , & Iv(0) = Iv_0
\end{eqnarray}
\begin{prop}
\label {Proposition 1}
For any $(u_0,v_0) \in H^s({\bf R}) \times H^s({\bf R})$ and $s \ge 1/3$ there 
exists $\delta \le 1$ and 
$\delta \sim (\|Iu_0\|_{H^1} + \|Iv_0\|_{H^1})^{-4-}$ , if $\beta \neq 0$ , and 
$\delta \sim 
(\|Iu_0\|_{H^1} + \|Iv_0\|_{H^1})^{-3-},$ if $\beta = 0$ , such that system 
(\ref{1.1}),(\ref{1.2}),(\ref{1.3}) has a unique local solution in the time 
interval $[0,\delta]$ with the 
property (dropping from now on 
$[0,\delta]$ from the notation):
$$ \|Iu\|_{X^{1,b}} + \|Iv\|_{Y^{1,b}} \le c (\|Iu_0\|_{H^1} + \|Iv_0\|_{H^1}) 
\, ,  $$ 
where $ b = \frac{1}{2} + $ . 
\end{prop}
{\bf Proof:} We want to construct a fixed point of the mapping 
$\tilde{S}=(\tilde{S}_0,\tilde{S}_1)$ 
induced by the integral equations belonging to the system 
(\ref{1.1''}),(\ref{1.2''}),(\ref{1.3''}):
\begin{eqnarray*}
\tilde{S}_0(Iu,Iv) & := & e^{it\partial_x^2} Iu_0 - i \int_0^t 
e^{i(t-s)\partial_x^2}[\alpha I(u(s)v(s)) + 
\beta I(|u(s)|^2 
u(s))] \, ds \, , \\
\tilde{S}_1(Iu,Iv) & := & e^{-t\partial_x^3} Iv_0 - \int_0^t 
e^{-(t-s)\partial_x^3} [I(v(s)\partial_x v(s)) 
- \gamma 
I(\partial_x(|u(s)|^2))] \, ds \, .
\end{eqnarray*}
The estimates for the nonlinearities in the previous proof carry over to 
corresponding 
estimates including the I-operator by the interpolation lemma of \cite{CKSTT5}, 
namely:
\begin{eqnarray*}
\|I(uv)\|_{X^{1,0}} & \le & c \|Iu\|_{X^{1,b}} \|Iv\|_{Y^{1,b}} \, ,  \\
\|I(|u|^2u)\|_{X^{1,0}} & \le & c \|Iu\|_{X^{1,b}}^3 \, , \\
\|I(v\partial_x v)\|_{Y^{1,-b'}} & \le & c \|Iv\|_{Y^{1,b}}^2 \quad  {\mbox for} 
\; b' \ge 
\max(\frac{1}{4}-\frac{s}{3},0) \; , \;\mbox{especially for} \; b' \ge 
\frac{1}{6} \, ,  \\
\|I\partial_x(|u|^2)\|_{Y^{1,-b'}} & \le & c \|Iu\|_{X^{1,b}}^2 \quad \mbox{for} 
\; b' > 
\max(\frac{1}{6},\frac{1}{2}-s) = \frac{1}{6} \, , \, \mbox{if} \; s \ge 
\frac{1}{3} \, .
\end{eqnarray*}
This implies
\begin{eqnarray*}
\|\tilde{S}_0(Iu,Iv)\|_{X^{1,b}} & \le & c \|Iu_0\|_{H^1} + c \alpha 
\|Iu\|_{X^{1,b}} \|Iv\|_{Y^{1,b}} 
\delta^{\frac{1}{2}-} +c\beta \|Iu\|_{X^{1,b}}^3 \delta^{\frac{1}{2}-} \, , \\
\|\tilde{S}_1(Iu,Iv)\|_{Y^{1,b}} & \le & c\|Iv_0\|_{H^1} + c \|Iv\|_{Y^{1,b}}^2 
\delta^{\frac{1}{3}-} 
+c\gamma 
\|Iu\|_{X^{1,b}}^2 \delta^{\frac{1}{3}-} \, .
\end{eqnarray*}
This gives the desired bounds, provided $ c\beta 
\delta^{\frac{1}{2}-}(\|Iu_0\|_{H^1}+\|Iv_0\|_{H^1})^2 < 1 
$ and $ c\delta^{\frac{1}{3}-}(\|Iu_0\|_{H^1} + \|Iv_0\|_{H^1}) < 1 $ . Thus our 
claimed choice of $\delta$ 
is possible.
\section{Conserved and almost conserved quantities}
Our system (\ref{1.1}),(\ref{1.2}),(\ref{1.3}) has the following conserved 
quantities (cf. \cite{Ts}):
\begin{eqnarray}
\label{3.0}
\hspace{-2em} M & := & \|u\| \, ,\\
\label{3.0a}
\hspace{-2em} L(u,v) & := & \alpha \|v\|^2 + 2\gamma \int Im(u\overline{u_x}) \, 
dx \, , \\
\label{3.0b}
 \hspace{-2em} E(u,v) & := & \alpha \gamma \int v|u|^2 dx + \gamma \|u_x\|^2 
+\frac{\alpha}{2}\|v_x\|^2 - 
\frac{\alpha}{6} 
\int v^3 dx + \frac{\beta \gamma}{2} \int |u|^4 dx \, .
\end{eqnarray}
Assume from now on $ \alpha \gamma > 0 $ .\\
These conservation laws imply a-priori bounds of the $H^1$ - norm of the 
solutions $u$ and $v$ as follows:
concerning $L$ we immediately get
\begin{equation}
\label{3.1}
|L(u,v)| \le c (\|v\|^2 + M \|u_x\|)
\end{equation}
and
\begin{equation}
\label{3.2}
\|v\|^2 \le c(|L|+M\|u_x\|) \, .
\end{equation}
Concerning $E$ we have by Gagliardo-Nirenberg
\begin{eqnarray*}
\int |v|^3 \, dx & \le & c \|v\|^{\frac{5}{2}} \|v_x\|^{\frac{1}{2}} \, \le \, 
c(|L|^{\frac{5}{4}} 
\|v_x\|^{\frac{1}{2}} +M\|u_x\|^{\frac{5}{4}} \|v_x\|^{\frac{1}{2}}) \\
& \le & c|L|^{\frac{5}{3}} + \epsilon \|v_x\|^2 +cM^{\frac{4}{3}} 
\|u_x\|^{\frac{5}{3}} \, \le \, 
c|L|^{\frac{5}{3}} + \epsilon (\|v_x\|^2 + \|u_x\|^2) + cM^8
\end{eqnarray*}
and
$$ \int |u|^4 \, dx \le \|u\|^3 \|u_x\| = M^3 \|u_x\| \le \epsilon \|u_x\|^2 + 
cM^6 $$
as well as
\begin{eqnarray*}
\int |vu^2| \, dx & \le & \|v\| \|u\|^{\frac{3}{2}} \|u_x\|^{\frac{1}{2}} \, \le 
\, c(|L|^{\frac{1}{2}} 
M^{\frac{3}{2}} \|u_x\|^{\frac{1}{2}} + M^2 \|u_x\|) \\
& \le & c(|L|^{\frac{2}{3}}M^2 + M^4) + \epsilon \|u_x\|^2 \, \le \, 
c(|L|^{\frac{5}{3}} + M^{\frac{10}{3}} 
+ M^4) + \epsilon \|u_x\|^2 \, .
\end{eqnarray*} 
This implies
$$ |\gamma| \|u_x\|^2 + \frac{|\alpha|}{2} \|v_x\|^2 \le |E| + 
c(|L|^{\frac{5}{3}} + M^8 +1) + \epsilon 
(\|u_x\|^2 + \|v_x\|^2) $$
and therefore
\begin{equation}
\label{3.2a}
\|u_x\|^2 + \|v_x\|^2 \le c(|E| + |L|^{\frac{5}{3}} + M^8 +1) \, .
\end{equation}
Similarly we also get
\begin{equation}
\label{3.3}
|E| \le c(\|u_x\|^2 + \|v_x\|^2 + |L|^{\frac{5}{3}} + M^8 +1) \, .
\end{equation}
The bounds (\ref{3.1}) and (\ref{3.3}) imply
\begin{eqnarray}
\nonumber
|E| & \le & c(\|u_x\|^2 + \|v_x\|^2 + \|v\|^{\frac{10}{3}} + M^{\frac{5}{3}} 
\|u_x\|^{\frac{5}{3}} + M^8 
+1)\\
\label{3.4}
& \le & c(\|u_x\|^2 + \|v_x\|^2 + \|v\|^{\frac{10}{3}} + M^{10} + 1) \, .
\end{eqnarray}
Finally, by (\ref{3.2}) and (\ref{3.2a}) we also have
\begin{equation}
\label{3.5}
\|v\|^2 \le c(|L|+M(|E|^{\frac{1}{2}} + |L|^{\frac{5}{6}} + M^4 +1)) \le c(|L| + 
M|E|^{\frac{1}{2}}+M^6+1)
\end{equation}
and thus
\begin{equation}
\label{3.6}
\|u\|_{H^1}^2 + \|v\|_{H^1}^2 \le c (|E| + |L|^{\frac{5}{3}} + M^8 +1) \, .
\end{equation}
These estimates imply an a-priori bound for the $H^1$ - norms of $u$ and $v$ for 
any data with finite 
energy $E$, finite $L$ and finite $\|u_0\|$. This is the case for $H^1$ - data 
$u_0$ and $v_0$. Thus we 
have from our local result (Theorem \ref{Theorem 1}):
\begin{theorem}
For data $(u_0,v_0) \in H^1({\bf R}) \times H^1({\bf R})$ and $\alpha \gamma > 
0$ there exists $b> 1/2$ 
such that (\ref{1.1}),(\ref{1.2}),(\ref{1.3}) has a unique global solution 
$(u,v) \in X^{1,b} \times 
Y^{1,b}$ with $(u,v) \in C^0({\bf R^+},H^1({\bf R}) \times H^1({\bf R}))$ .
\end{theorem}

A crucial role is played by the modified functionals $E(Iu,Iv)$ and $L(Iu,Iv)$, 
which are "almost" 
conserved, i.e. their growth in time is controllable. Using the modified system 
(\ref{1.1''}),(\ref{1.2''}),(\ref{1.3''}), an elementary calculation shows
\begin{eqnarray}
\nonumber
\lefteqn{ \frac{d}{dt} E(Iu,Iv) } \\
\nonumber
& = & \alpha \left[\int(I(vv_x)-IvIv_x)Iv_{xx} dx + \frac{1}{2}\int 
(Iv)^2(I(vv_x)-IvIv_x)dx \right] \\
\nonumber
& & +2\beta\gamma Im\int(I(|u|^2u)_x - ((Iu)^2I\overline{u})_x)I\overline{u}_x 
dx \\
\nonumber
& & +\alpha\gamma [\int |Iu|^2(Iv Iv_x - I(vv_x))dx + \int (|Iu|^2-I(|u|^2))Iv 
Iv_x \, dx \\
& & \hspace{0.5cm} + \int Iv_{xx}(|Iu|^2_x - I(|u|^2)_x)dx 
  - 2 Im \int Iu_x(I(\overline{u}v)_x - (I\overline{u}Iv)_x)dx ] 
  \label{3.1'} \\
\nonumber
& & + \alpha \gamma^2 \int (I(|u|^2)_x - (|Iu|^2)_x)|Iu|^2 dx +2\alpha^2 \gamma 
Im \int 
(I(\overline{u}v)-I\overline{u}Iv)Iv Iu dx \\
\nonumber
& & +2\beta^2\gamma Im \int Iu 
(I\overline{u})^2(I(|u|^2u)-(Iu)^2I\overline{u})dx \\
\nonumber
& & -2\alpha\beta\gamma[Im \hspace{-0.2em}\int\hspace{-0.3em} Iv Iu 
(I(|u|^2\overline{u})-Iu(I\overline{u})^2)dx + Im\hspace{-0.2em} 
\int\hspace{-0.2em} 
(Iu)^2I\overline{u}(I(\overline{u}v)-I\overline{u}Iv)dx] \\
& =: & \sum_{j=1}^{12} I_j \, .
\nonumber
\end{eqnarray}
Similarly
\begin{eqnarray}
\nonumber
\lefteqn{\frac{d}{dt} L(Iu,Iv)
 =  2\alpha \int Iv(Iv Iv_x - I(vv_x)) dx} \\
\nonumber & & +2\alpha\gamma\left[\int Iv(I(|u|^2)_x - (|Iu|^2)_x)dx + 2 Re\int 
I\overline{u}_x(IuIv-I(uv))dx\right] \\
\label{3.2'}
& & + 2\beta\gamma Re \int ((Iu)^2I\overline{u} - I(u^2 
\overline{u}))I\overline{u}_x dx \\
\nonumber
& & =: \sum_{j=1}^4 J_j \, .
\end{eqnarray}
\section{Estimates for the modified energy functional}
We need exact control of the increment of the modified energy.
\begin{prop}
\label{Proposition 4.1}
If $(u,v)$ is a solution of (\ref{1.1}),(\ref{1.2}),(\ref{1.3}) on $[0,\delta]$ 
in the sense of Proposition 
\ref{Proposition 1}, then the following estimate holds for $N \ge 1$ and $ s > 
1/2 $ :
\begin{eqnarray*}
\lefteqn{|E(Iu(\delta),Iv(\delta))-E(Iu(0),Iv(0))|} \\
& \hspace{-0.8em} \le & \hspace{-0.8em} c[(N^{-1+} \delta^{\frac{1}{2}-} + 
N^{-\frac{7}{4}+})(\|Iu\|_{X^{1,\frac{1}{2}+}}^3 + 
\|Iv\|_{Y^{1,\frac{1}{2}+}}^3) \\
& & \hspace{-1em} + N^{-2+}(\|Iu\|^4_{X^{1,\frac{1}{2}+}} + 
\|Iv\|^4_{Y^{1,\frac{1}{2}+}}) + 
N^{-3+}\|Iu\|^4_{X^{1,\frac{1}{2}+}}(\|Iu\|_{X^{1,\frac{1}{2}+}}^2 + 
\|Iv\|_{Y^{1,\frac{1}{2}+}})] \, .
\end{eqnarray*}
\end{prop}
{\bf Proof:} Integrating (\ref{3.1'}) over $t \in [0,\delta]$ we have to 
estimate the various terms on the 
r.h.s. Here and in the sequel we assume w.l.o.g. the Fourier transforms of all 
the functions to be 
nonnegative and ignore the appearance of complex conjugates if this is 
irrelevant for the argument. We use 
dyadic decompositions w.r. to the frequencies $|\xi_j| \sim N_j = 2^k $ 
$(k=0,1,2,...)$ in many places, so 
that we need extra factors $N_j^{0-}$ everywhere in order to sum the dyadic 
pieces.\\
\underline{Estimate of $I_1$:} We have to show
\begin{eqnarray}
\nonumber
& & \hspace{-2em}\int_0^{\delta} \int_* \left| 
\frac{m(\xi_1+\xi_2)-m(\xi_1)m(\xi_2)}{m(\xi_1)m(\xi_2)}\right| |\xi_1+\xi_2| 
\widehat{v_1}(\xi_1,t) 
|\xi_2| \widehat{v_2}(\xi_2,t) |\xi_3| \widehat{v_3}(\xi_3,t) d\xi dt \\
\label{4.1}
& & \le  c(N^{-1+} \delta^{\frac{1}{2}-} + N^{-\frac{7}{4}+}) \prod_{i=1}^3 
\|v_i\|_{Y^{1,\frac{1}{2}+}} \, 
.
\end{eqnarray}
Here and in the sequel * denotes integration over the set $\sum \xi_i = 0$. 
Typically, because of the 
multiplier term we can assume $|\xi_1| \ge N$ or $|\xi_2| \ge N$, and because of 
the convolution 
constraint the two largest frequencies are equivalent.\\
\underline{Case 1:} $|\xi_1| << |\xi_2| \, , \, |\xi_1| \le cN \, , \, |\xi_2| 
\ge cN \, , \, |\xi_3| \sim 
|\xi_2| $ .\\
The multiplier term is estimated using the mean value theorem by
$$ \left| \frac{m(\xi_1+\xi_2)-m(\xi_2)}{m(\xi_2)}\right| \le c \frac{|(\nabla 
m)(\xi_2) 
\xi_1|}{|m(\xi_2)|} \le c \frac{N_1}{N_2} \, . $$
Thus the integral is bounded by use of (\ref{g}) and (\ref{0}):
\begin{eqnarray*}
\lefteqn{ c \frac{N_1}{N_2} \int_0^{\delta} \int_* |\xi_1+\xi_2|^{\frac{1}{2}} 
|\xi_1-\xi_2|^{\frac{1}{2}} 
\widehat{v_1}(\xi_1,t) |\xi_2| \widehat{v_2}(\xi_2,t) |\xi_3| 
\widehat{v_3}(\xi_3,t) d\xi dt} \\
& &\hspace{-1em}\le  c \frac{N_1}{N_2} \|D_x^{\frac{1}{2}}I_-^{\frac{1}{2}}(v_1, 
D_x v_2)\|_{L^2_{xt}} 
\|D_x 
v_3\|_{L^2_{xt}} 
 \le c \frac{N_1}{N_2} \|v_1\|_{Y^{0,\frac{1}{2}+}} \|D_x 
v_2\|_{Y^{0,\frac{1}{2}+}} \|D_x v_3\|_{L^2_{xt}} 
\\
& & \hspace{-1em}\le c \frac{N_1}{N_2} \frac{1}{N_1} \delta^{\frac{1}{2}-} 
\prod_{i=1}^3 
\|v_i\|_{Y^{1,\frac{1}{2}+}}  \le c N^{-1+} N_{max}^{0-} \delta^{\frac{1}{2}-} 
\prod_{i=1}^3 
\|v_i\|_{Y^{1,\frac{1}{2}+}} \, . 
\end{eqnarray*}
\underline{Case 2:} $|\xi_1| << |\xi_2| \, , \, |\xi_1|,|\xi_2| \ge N$ .\\
Here we avoid dyadic decompositions and estimate the multiplier by \\ 
$\frac{c}{m(\xi_1)} \le 
c(\frac{|\xi_1|}{N})^{\frac{1}{2}}$ . Similarly as in case 1 we control the 
integral by
\begin{eqnarray*}
\lefteqn{cN^{-\frac{1}{2}}\int_0^{\delta}\int_* |\xi_1+\xi_2|^{\frac{1}{2}} 
|\xi_1-\xi_2|^{\frac{1}{2}} 
|\xi_1|^{\frac{1}{2}} \widehat{v_1}(\xi_1,t) |\xi_2| \widehat{v_2}(\xi_2,t) 
|\xi_3| \widehat{v_3}(\xi_3,t) 
d\xi dt} \\
& & \le  cN^{-\frac{1}{2}} \|D_x^{\frac{1}{2}}v_1\|_{Y^{0,\frac{1}{2}+}} \|D_x 
v_2\|_{Y^{0,\frac{1}{2}+}} 
\|D_x v_3\|_{L^2_{xt}} \le cN^{-1}\delta^{\frac{1}{2}-} \prod_{i=1}^3 
\|v_i\|_{Y^{1,\frac{1}{2}+}} \, .
\end{eqnarray*}
Similarly we treat the case $|\xi_1| >> |\xi_2|$ , so that it remains to 
consider\\
\underline{Case 3:} $|\xi_1| \sim |\xi_2| \ge N \, .$  \\
The multiplier is bounded by $cN_1 N^{-1}$ and thus we get using (\ref{f}):
\begin{eqnarray*}
\lefteqn{c N_1 N^{-1} \|D_x(v_1 D_x v_2)\|_{Y^{0,-\frac{1}{2}-}} \|D_x 
v_3\|_{Y^{0,\frac{1}{2}+}}} \\ & 
&\le c N_1 N^{-1} \|v_1\|_{Y^{-\frac{3}{8}+,\frac{1}{2}+}} \|D_x 
v_2\|_{Y^{-\frac{3}{8}+,\frac{1}{2}+}} 
\|D_x v_3\|_{Y^{0,\frac{1}{2}+}} \\
& & \le c N_1 N^{-1} N_1^{-1-\frac{3}{8}+} N_2^{-\frac{3}{8}+} \prod_{i=1}^3 
\|v_i\|_{Y^{1,\frac{1}{2}+}} 
\le c 
N^{-\frac{7}{4}+} N_{max}^{0-} \prod_{i=1}^3 \|v_i\|_{Y^{1,\frac{1}{2}+}} \, .
\end{eqnarray*}
\underline{Estimate of $I_2$:} It is sufficient to show
\begin{eqnarray*}
& & \hspace{-2em}
\int_0^{\delta} \int_* 
\left|\frac{m(\xi_3+\xi_4)-m(\xi_3)m(\xi_4)}{m(\xi_3)m(\xi_4)}\right| 
\widehat{v_1}(\xi_1,t) \widehat{v_2}(\xi_2,t) \widehat{v_3}(\xi_3,t) |\xi_4| 
\widehat{v_4}(\xi_4,t) \, d\xi 
dt \\
& & \le c N^{-2+} \prod_{i=1}^4 \|v_i\|_{Y^{1,\frac{1}{2}+}} \, .
\end{eqnarray*}
At least two of the $N_i$ are $ \ge cN$. Assume w.l.o.g. $N_1 \ge N_2 \ge N_3$ 
such that $N_1 \ge cN$ . 
Then we get by use of (\ref{e}) the bound
\begin{eqnarray*}
\lefteqn{ c N_{max} N^{-1} \|v_1 v_2\|_{L^2_{xt}} \|v_3 D_x v_4\|_{L^2_{xt}} } 
\\
& \le & c N_{max} N^{-1} \|v_1\|_{Y^{-\frac{1}{2},\frac{1}{2}+}} 
\|v_2\|_{Y^{\frac{1}{4},\frac{1}{2}+}}\|v_3\|_{Y^{\frac{1}{4},\frac{1}{2}+}}\|D_
x 
v_4\|_{Y^{-\frac{1}{2},\frac{1}{2}+}} \\
& \le & c N_{max} N^{-1} N_1^{-\frac{3}{2}} \langle N_2 \rangle ^{-\frac{3}{4}} 
\langle N_3 \rangle 
^{-\frac{3}{4}} \langle N_4 \rangle ^{-\frac{1}{2}} \prod_{i=1}^4 
\|v_i\|_{Y^{1,\frac{1}{2}+}} \\
& \le & c N^{-2+} N_{max}^{0-} \prod_{i=1}^4 \|v_i\|_{Y^{1,\frac{1}{2}+}} \, .
\end{eqnarray*}
\underline{Estimate of $I_3$:} We have to show
\begin{eqnarray*}
& & \hspace{-2em}
\int_0^{\delta} \int_* \left| 
\frac{m(\xi_1+\xi_2+\xi_3)-m(\xi_1)m(\xi_2)m(\xi_3)}{m(\xi_1)m(\xi_2)m(\xi_3)}
\right| \times \\ & & 
\hspace{1em} \times \widehat{u_1}(\xi_1,t) \widehat{u_2}(\xi_2,t) 
\widehat{u_3}(\xi_3,t) \xi_4^2 
\widehat{u_4}(\xi_4,t)  d\xi dt \,
 \le  c N^{-2+} \prod_{i=1}^4 \|u_i\|_{X^{1,\frac{1}{2}+}} \, .
\end{eqnarray*}
\underline{Case 1:} $ |\xi_1| \sim |\xi_2| \sim |\xi_3| \sim|\xi_4| \ge cN \, 
.$\\
We get the bound by use of Strichartz:
\begin{eqnarray*}
\lefteqn{ c (N_1/N)^{\frac{3}{2}} \|u_1\|_{L^4_{xt}} \|u_2\|_{L^4_{xt}} 
\|u_3\|_{L^4_{xt}} \|D_x^2 
u_4\|_{L^4_{xt}}} \\
& \le & c (N_1/N)^{\frac{3}{2}} (N_1 N_2 N_3)^{-1} N_4 
\prod_{i=1}^4{\|u_i\|_{X^{1,\frac{1}{2}+}}} \\
& \le & c N^{-2+} N_{max}^{0-} \prod_{i=1}^4{\|u_i\|_{X^{1,\frac{1}{2}+}}} \, .
\end{eqnarray*}
\underline{Case 2:} Two of the frequencies are $\ge cN$, the most difficult case 
is $|\xi_4| \ge cN$ and, 
say, $|\xi_1| \sim |\xi_4|$ , $ |\xi_2|,|\xi_3| << |\xi_1|,|\xi_4|$ . \\
a. $ |\xi_2| \ge N $ (or similarly $|\xi_3| \ge N$). \\
In this case the multiplier is bounded by $c\langle (N_3/N)^{\frac{1}{2}} 
\rangle (N_2/N)^{\frac{1}{2}} $ , 
and using (\ref{c}) we get the bound
\begin{eqnarray*}
\lefteqn{ c \langle (N_3/N)^{\frac{1}{2}} \rangle  (N_2/N)^{\frac{1}{2}} \|u_1 
u_2\|_{L^2_{xt}} \|u_3 D_x^2 
u_4\|_{L^2_{xt}} } \\
& \le & c \langle (N_3/N)^{\frac{1}{2}} \rangle (N_2/N)^{\frac{1}{2}} 
\|D_x^{-\frac{1}{2}} 
u_1\|_{X^{0,\frac{1}{2}+}} \|u_2\|_{X^{0,\frac{1}{2}+}} 
\|u_3\|_{X^{0,\frac{1}{2}+}} 
\|D_x^{\frac{3}{2}}u_4\|_{X^{0,\frac{1}{2}+}} \\
& \le & c \langle (N_3/N)^{\frac{1}{2}} \rangle (N_2/N)^{\frac{1}{2}} 
N_1^{-\frac{3}{2}} N_2^{-1} \langle 
N_3\rangle ^{-1} N_4^{\frac{1}{2}} \prod_{i=1}^4 \|u_i\|_{X^{1,\frac{1}{2}+}} \\
& \le & c N^{-2+} N_{max}^{0-} \prod_{i=1}^4 \|u_i\|_{X^{1,\frac{1}{2}+}} \, .
\end{eqnarray*}
b. $|\xi_2|,|\xi_3| \le N. $ \\
By use of the mean value theorem the multiplier is bounded by
$$ \left| \frac{m(\xi_1+\xi_2+\xi_3)-m(\xi_1)}{m(\xi_1)} \right| = \left| 
\frac{(\nabla 
m)(\xi_1)\cdot(\xi_2 + \xi_3)}{m(\xi_1)} \right| \le \frac{N_2 + N_3}{N_1} $$
and exactly as in case a the claimed estimate follows.\\
\underline{Case 3:} Three of the frequencies are equivalent, say $|\xi_1| \sim 
|\xi_2| \sim |\xi_4| \ge cN 
$ , $ |\xi_3| << |\xi_1|,|\xi_2|,|\xi_4|$ . Then two  of the large frequencies 
have different sign, say 
$\xi_1$ and $\xi_2$, so that $|\xi_1 - \xi_2|^{\frac{1}{2}} \sim 
|\xi_1|^{\frac{1}{2}}$ (the other cases 
are treated similarly). We get the bound, using (\ref{h}) and (\ref{c}):
\begin{eqnarray*}
\lefteqn{ c (N_1/N)^{\frac{1}{2}} (N_2/N)^{\frac{1}{2}} \langle 
(N_3/N)^{\frac{1}{2}} \rangle \|u_1 
u_2\|_{L^2_{xt}} \|u_3 D_x^2 u_4\|_{L^2_{xt}} } \\
& \le & c (N_1/N)^{\frac{1}{2}} (N_2/N)^{\frac{1}{2}} \langle 
(N_3/N)^{\frac{1}{2}} \rangle 
\|I_-^{\frac{1}{2}}(D_x^{-\frac{1}{2}} u_1, u_2)\|_{L^2_{xt}} \|u_3 D_x^2 
u_4\|_{L^2_{xt}} \\ 
& \le & c (N_1/N)^{\frac{1}{2}} (N_2/N)^{\frac{1}{2}} \langle 
(N_3/N)^{\frac{1}{2}} \rangle \times \\
& &  \hspace{10em}\times \, \|D_x^{-\frac{1}{2}}u_1\|_{X^{0,\frac{1}{2}+}} 
\|u_2\|_{X^{0,\frac{1}{2}+}} 
\|u_3\|_{X^{0,\frac{1}{2}+}} \|D_x^{\frac{3}{2}}u_1\|_{X^{0,\frac{1}{2}+}} \\
& \le & c (N_1/N)^{\frac{1}{2}} (N_2/N)^{\frac{1}{2}} \langle 
(N_3/N)^{\frac{1}{2}} \rangle 
N_1^{-\frac{3}{2}} N_2^{-1} \langle N_3 \rangle ^{-1} N_4^{\frac{1}{2}} 
\prod_{i=1}^4 
\|u_i\|_{X^{1,\frac{1}{2}+}} \\
& \le & c N^{-2+} N_{max}^{0-} \prod_{i=1}^4 \|u_i\|_{X^{1,\frac{1}{2}+}} \, . 
\end{eqnarray*}
\underline{Estimate of $I_4$:} We have to show
\begin{eqnarray*}
& &  \int_0^{\delta} \int_* \left| 
\frac{m(\xi_1+\xi_2)-m(\xi_1)m(\xi_2)}{m(\xi_1)m(\xi_2)} \right| 
\widehat{v_1}(\xi_1,t) |\xi_2| \widehat{v_2}(\xi_2,t) \widehat{u_3}(\xi_3,t) 
\widehat{u_4}(\xi_4,t) d\xi dt 
 \\
& & \quad \le c N^{-2+} \|u_1\|_{X^{1,\frac{1}{2}+}} 
\|u_2\|_{X^{1,\frac{1}{2}+}} 
\|v_1\|_{Y^{1,\frac{1}{2}+}} \|v_2\|_{Y^{1,\frac{1}{2}+}} \, .
\end{eqnarray*}
We bound the multiplier by $c N_{max}N^{-1}$. We have $|\xi_1|\ge N$ or $|\xi_2| 
\ge N$. Assume the more 
difficult case $|\xi_2| \ge N$.\\
\underline{Case 1:} Exactly two of the $N_i$ are $\ge cN$, thus w.l.o.g. $N_3 << 
N_2$.\\
Using Lemma \ref{Lemma 3} and Strichartz we get the bound
\begin{eqnarray*}
\lefteqn{ c N_{max} N^{-1} \|(D_x v_2) u_3\|_{L^2_{xt}} \|v_1\|_{L^4_{xt}} 
\|u_4\|_{L^4_{xt}} } \\
& \le & c N_{max} N^{-1} \|v_2\|_{Y^{0,\frac{1}{2}+}} 
\|u_3\|_{X^{0,\frac{1}{2}+}} 
\|v_1\|_{Y^{0,\frac{1}{2}+}} \|u_4\|_{X^{0,\frac{1}{2}+}} \\
& \le & c N_{max} N^{-1} (N_1 N_2 N_3 N_4)^{-1} \prod_{i=3}^4 
\|u_i\|_{X^{1,\frac{1}{2}+}} \prod_{i=1}^2 
\|v_i\|_{Y^{1,\frac{1}{2}+}} \\
& \le & c N^{-2+} N_{max}^{0-} \prod_{i=3}^4 \|u_i\|_{X^{1,\frac{1}{2}+}} 
\prod_{i=1}^2 
\|v_i\|_{Y^{1,\frac{1}{2}+}} \, .
\end{eqnarray*}
\underline{Case 2:} At least three of the $N_i$ are $\ge cN$.\\
In this case Strichartz directly gives the bound
\begin{eqnarray*}
\lefteqn{c (\frac{N_{\max}}{N})^{1-\epsilon} \|v_1\|_{L^4_{xt}}  \|D_x 
v_2\|_{L^4_{xt}} \|u_3\|_{L^4_{xt}} 
\|u_4\|_{L^4_{xt}} } \\
& \le & c (\frac{N_{\max}}{N})^{1-\epsilon} (N_1N_3N_4)^{-1} \prod_{i=3}^4 
\|u_i\|_{X^{1,\frac{1}{2}+}} 
\prod_{i=1}^2 \|v_i\|_{Y^{1,\frac{1}{2}+}} \\
& \le & c N^{-2+} N_{max}^{0-} \prod_{i=3}^4 \|u_i\|_{X^{1,\frac{1}{2}+}} 
\prod_{i=1}^2 
\|v_i\|_{Y^{1,\frac{1}{2}+}} \, . 
\end{eqnarray*}
\underline{Estimate of $I_5$:} We have to show
\begin{eqnarray*}
& &  \int_0^{\delta} \int_* \left| 
\frac{m(\xi_1+\xi_2)-m(\xi_1)m(\xi_2)}{m(\xi_1)m(\xi_2)} \right| 
\widehat{u_1}(\xi_1,t)  \widehat{u_2}(\xi_2,t) \widehat{v_3}(\xi_3,t) |\xi_4| 
\widehat{v_4}(\xi_4,t) d\xi dt 
 \\
& & \quad \le c N^{-2+} \|u_1\|_{X^{1,\frac{1}{2}+}} 
\|u_2\|_{X^{1,\frac{1}{2}+}} 
\|v_1\|_{Y^{1,\frac{1}{2}+}} \|v_2\|_{Y^{1,\frac{1}{2}+}} \, .
\end{eqnarray*}
This can be shown similarly as the previous case.\\   
\underline{Estimate of $I_6$:} We want to show
\begin{eqnarray*}
& & \int_0^{\delta} \int_* 
\left|\frac{m(\xi_1+\xi_2)-m(\xi_1)m(\xi_2)}{m(\xi_1)m(\xi_2)} \right| 
|\xi_1+\xi_2| \widehat{u_1}(\xi_1,t)\widehat{\overline{u_2}}(\xi_2,t) \xi_3^2 
\widehat{v_3}(\xi_3,t) d\xi 
dt \\
& & \quad \le c N^{-1+} \delta^{\frac{1}{2}-} \|u_1\|_{X^{1,\frac{1}{2}+}}  
\|u_2\|_{X^{1,\frac{1}{2}+}}  
\|v_3\|_{Y^{1,\frac{1}{2}+}} \, .
\end{eqnarray*}
\underline{Case 1:} $|\xi_3| \sim |\xi_1| \ge cN $ , $|\xi_2| << |\xi_1|,|\xi_3| 
$ (thus $|\xi_1+\xi_2| 
\sim |\xi_1|$). \\
a. $|\xi_2| \le N$ \\
The multiplier is bounded by $c|\frac{(\nabla m)(\xi_1)\xi_2}{m(\xi_1)}| \le c 
\frac{N_2}{N_1}$ . This 
implies by Lemma \ref{Lemma 3} the bound 
\begin{eqnarray*}
\lefteqn{ c N_2 N_1^{-1} \|D_x u_1\|_{L^2_{xt}} \|u_2 D_x^2 v_3\|_{L^2_{xt}} } 
\\
& \le & c N_2 N_1^{-1} \|D_x u_1\|_{L^2_{xt}} \|u_2\|_{X^{0,\frac{1}{2}+}} \|D_x 
v_3\|_{Y^{0,\frac{1}{2}+}} 
\\
& \le & c N_2 N_1^{-1} \delta^{\frac{1}{2}-} \|u_1\|_{X^{1,\frac{1}{2}+}} 
N_2^{-1} 
\|u_2\|_{X^{1,\frac{1}{2}+}} \|v_3\|_{Y^{1,\frac{1}{2}+}} \\
& \le & c N^{-1+} N_{max}^{0-} \delta^{\frac{1}{2}-} 
\|u_1\|_{X^{1,\frac{1}{2}+}} 
\|u_2\|_{X^{1,\frac{1}{2}+}} \|v_3\|_{Y^{1,\frac{1}{2}+}} \, .
\end{eqnarray*}
b. $|\xi_2| \ge N$ \\
In this case we perform no dyadic decomposition and estimate the multiplier by 
$c(\frac{|\xi_2|}{N})^{\frac{1}{2}}$. This implies the following bound for the 
integral by Lemma \ref{Lemma 
3}:
\begin{eqnarray*}
\lefteqn{ cN^{-\frac{1}{2}} \int_0^{\delta} \int_* |\xi_1| 
\widehat{u_1}(\xi_1,t) |\xi_2|^{\frac{1}{2}} 
\widehat{u_2}(\xi_2,t) \xi_3^2 \widehat{v_3}(\xi_3,t) d\xi dt } \\
& \le & c N^{-\frac{1}{2}} \|D_x u_1\|_{L^2_{xt}} \|D_x^{\frac{1}{2}} u_2 D_x^2 
v_3\|_{L^2_{xt}} \\
& \le & c N^{-\frac{1}{2}} \|D_x u_1\|_{L^2_{xt}} \|D_x^{\frac{1}{2}} 
u_2\|_{X^{0,\frac{1}{2}+}} \| D_x 
v_3\|_{Y^{0,\frac{1}{2}+}} \\
& \le & c N^{-\frac{1}{2}} \delta^{\frac{1}{2}-} N^{-\frac{1}{2}}  
\|u_1\|_{X^{1,\frac{1}{2}+}} 
\|u_2\|_{X^{1,\frac{1}{2}+}} \|v_3\|_{Y^{1,\frac{1}{2}+}} \, .
\end{eqnarray*}
\underline{Case 2:} $|\xi_1| \sim |\xi_2| \ge cN $ , $|\xi_3|^2 >> |\xi_2| $ , 
thus $|\xi_2+\xi_1| \le 
c|\xi_1|$ .\\
We get the bound by Lemma \ref{Lemma 3}:
\begin{eqnarray*}
\lefteqn{ c (\frac{N_2}{N})^{1-\epsilon} \|D_x u_1\|_{L^2_{xt}} \|u_2 D_x^2 
v_3\|_{L^2_{xt}} } \\
& \le & c (\frac{N_2}{N})^{1-\epsilon} \|D_x u_1\|_{L^2_{xt}} 
\|u_2\|_{X^{0,\frac{1}{2}+}} \|D_x 
v_3\|_{Y^{0,\frac{1}{2}+}} \\
& \le & c (\frac{N_2}{N})^{1-\epsilon} \delta^{\frac{1}{2}-} N_2^{-1} 
\|u_1\|_{X^{1,\frac{1}{2}+}} 
\|u_2\|_{X^{1,\frac{1}{2}+}} \|v_3\|_{Y^{1,\frac{1}{2}+}} \\
& \le & c N^{-1+} N_{max}^{0-} \delta^{\frac{1}{2}-} 
\|u_1\|_{X^{1,\frac{1}{2}+}} 
\|u_2\|_{X^{1,\frac{1}{2}+}} \|v_3\|_{Y^{1,\frac{1}{2}+}} \, . 
\end{eqnarray*}
\underline{Case 3:} $|\xi_1| \sim|\xi_2| \ge cN $ , $ |\xi_3|^2 \le c|\xi_1| 
\sim c|\xi_2| $ . \\
Using $|\xi_1+\xi_2| \xi_3^2 \le c|\xi_1+\xi_2|^{\frac{1}{2}}|\xi_1||\xi_3|$ and 
the multiplier bound 
$c(\frac{N_2}{N})^{\frac{1}{2}}$ we estimate the integral by use of (\ref{d}):
\begin{eqnarray*}
\lefteqn{ c(\frac{N_2}{N})^{1-\epsilon} \|D_x^{\frac{1}{2}}(D_x u_1 
\overline{u_2})\|_{L^2_{xt}} \|D_x 
v_3\|_{L^2_{xt}} } \\
& \le & c(\frac{N_2}{N})^{1-\epsilon} \|D_x u_1\|_{X^{0,\frac{1}{2}+}} 
\|u_2\|_{X^{0,\frac{1}{2}+}} \|D_x 
v_3\|_{L^2_{xt}} \\
& \le & c(\frac{N_2}{N})^{1-\epsilon} N_2^{-1} \delta^{\frac{1}{2}-} 
\|u_1\|_{X^{1,\frac{1}{2}+}} 
\|u_2\|_{X^{1,\frac{1}{2}+}} \|v_3\|_{Y^{1,\frac{1}{2}+}} \\
& \le & c N^{-1+} N_{max}^{0-} \delta^{\frac{1}{2}-} 
\|u_1\|_{X^{1,\frac{1}{2}+}} 
\|u_2\|_{X^{1,\frac{1}{2}+}} \|v_3\|_{Y^{1,\frac{1}{2}+}} \, .
\end{eqnarray*}
\underline{Estimate of $I_7$:} We show
\begin{eqnarray*}
& & \int_0^{\delta} \int_* 
\left|\frac{m(\xi_1+\xi_2)-m(\xi_1)m(\xi_2)}{m(\xi_1)m(\xi_2)}\right||\xi_1+\xi_
2| \widehat{u_1}(\xi_1,t) 
\widehat{v_2}(\xi_2,t) |\xi_3| \widehat{u_3}(\xi_3,t) d\xi dt \\
& & \quad \le c N^{-1} \delta^{\frac{1}{2}-} \|u_1\|_{X^{1,\frac{1}{2}+}} 
\|v_2\|_{Y^{1,\frac{1}{2}+}} 
\|u_3\|_{X^{1,\frac{1}{2}+}} \, .
\end{eqnarray*}
\underline{Case 1:} $|\xi_1| >> |\xi_2| \ge cN $ . \\
Using no dyadic decomposition and the multiplier bound 
$c(\frac{|\xi_2|}{N})^{\frac{1}{2}}$ we get the 
bound by Strichartz:
\begin{eqnarray*}
& & c N^{-\frac{1}{2}} \int_0^{\delta} \int_* |\xi_1+\xi_2| 
\widehat{u_1}(\xi_1,t) |\xi_2|^{\frac{1}{2}} 
\widehat{v_2}(\xi_2,t) |\xi_3| \widehat{u_3}(\xi_3,t) d\xi dt \\
& & \quad \le c N^{-\frac{1}{2}} \|D_xu_1\|_{L^4_{xt}} 
\|D_x^{\frac{1}{2}}v_2\|_{L^4_{xt}} 
\|D_xu_3\|_{L^2_{xt}} \\
& & \quad \le c N^{-\frac{1}{2}} \|u_1\|_{X^{1,\frac{1}{2}+}} 
N^{-\frac{1}{2}}\|v_2\|_{Y^{1,\frac{1}{2}+}} 
\|u_3\|_{X^{1,\frac{1}{2}+}} \delta^{\frac{1}{2}-} \, .
\end{eqnarray*} 
\underline{Case 2:} $|\xi_1| \sim |\xi_2| \ge cN $ . \\
Again using no dyadic decomposition and the multiplier bound 
$c\frac{|\xi_2|}{N}$ we get similarly as in 
case 1 the bound
$$ cN^{-1} \|D_xu_1\|_{L^4_{xt}} \|D_xv_2\|_{L^4_{xt}} \|D_xu_3\|_{L^2_{xt}} \le 
c N^{-1} 
\delta^{\frac{1}{2}-} \|u_1\|_{X^{1,\frac{1}{2}+}} \|v_2\|_{Y^{1,\frac{1}{2}+}} 
\|u_3\|_{X^{1,\frac{1}{2}+}} \, . $$
\underline{Case 3:} $|\xi_1| \ge cN $ , $ |\xi_2| \le N $ , thus $N_1 \sim 
N_{max}$ . \\
We bound the multiplier by $c \frac{N_2}{N_1}$ and get by Strichartz an integral 
bound
$$c \frac{N_2}{N_1} \|D_xu_1\|_{L^4_{xt}} \|v_2\|_{L^4_{xt}} \|D_x 
u_3\|_{L^2_{xt}} \le c N_1^{-1} 
\delta^{\frac{1}{2}-} \|u_1\|_{X^{1,\frac{1}{2}+}} \|v_2\|_{Y^{1,\frac{1}{2}+}} 
\|u_3\|_{X^{1,\frac{1}{2}+}} \, . $$
The remaining cases are handled similarly by exchanging the roles of $\xi_1$ and 
$\xi_2$.\\
\underline{Estimate of $I_8$:} We want to show
$$  \int_0^{\delta} \int_* 
\left|\frac{m(\xi_1+\xi_2)-m(\xi_1)m(\xi_2)}{m(\xi_1)m(\xi_2)}\right| 
|\xi_1+\xi_2| \prod_{i=1}^4 \widehat{u_i}(\xi_i,t) \, d\xi dt \le  c N^{-2+} 
\prod_{i=1}^4 
\|u_i\|_{X^{1,\frac{1}{2}+}}. $$
\underline{Case 1:} At least three of the $|\xi_i|$ are $\ge cN$ , $|\xi_1| \ge 
|\xi_2|$ w.l.o.g. \\
Estimate the multiplier by $cN_{max}N^{-1}$ and use Strichartz to control the 
integral by
\begin{eqnarray*}
c N_{max} N^{-1} \prod_{i=2}^4 \|u_i\|_{L^4_{xt}} \|D_x u_1\|_{L^4_{xt}} & \le & 
c N_{max} N^{-1} 
\prod_{i=2}^4 \langle N_i \rangle ^{-1} \prod_{i=1}^4 
\|u_i\|_{X^{1,\frac{1}{2}+}} \\
& \le & c N^{-2+} N_{max}^{0-} \prod_{i=1}^4 \|u_i\|_{X^{1,\frac{1}{2}+}} \, .
\end{eqnarray*}
\underline{Case 2:} Exactly two of the $|\xi_i|$ are $\ge cN$, the others $ << N 
$ , e.g. $|\xi_1| \sim 
|\xi_2| \ge cN $ and $ |\xi_3|,|\xi_4| << N $ . \\
Estimate the multiplier by $cN_{max}N^{-1}$ and use (\ref{c}) to bound the 
integral by
\begin{eqnarray*}
\lefteqn{ cN_{max}N^{-1}\|(D_x u_1) u_3\|_{L^2_{xt}} \|u_2 u_4\|_{L^2_{xt}} } \\
 & \le & c N_{max}N^{-1} \|D_x^{\frac{1}{2}} u_1\|_{X^{0,\frac{1}{2}+}} 
\|u_3\|_{X^{0,\frac{1}{2}+}} 
\|D_x^{-\frac{1}{2}} u_2\|_{X^{0,\frac{1}{2}+}} \|u_4\|_{X^{0,\frac{1}{2}+}} \\ 
&
\le & c N_{max}N^{-1} N_1^{-\frac{1}{2}} N_2^{-\frac{3}{2}} \prod_{i=1}^4 
\|u_i\|_{X^{1,\frac{1}{2}+}} 
\quad \le \quad c N^{-2+} N_{max}^{0-} \prod_{i=1}^4 
\|u_i\|_{X^{1,\frac{1}{2}+}} .
\end{eqnarray*}
\underline{Estimate of $I_9$:} We control the integral
$$ \int_0^{\delta} \int_* 
\left|\frac{m(\xi_1+\xi_2)-m(\xi_1)m(\xi_2)}{m(\xi_1)m(\xi_2)}\right| 
\widehat{u_1}(\xi_1,t)
\widehat{v_2}(\xi_2,t) \widehat{u_3}(\xi_3,t) \widehat{v_4}(\xi_4,t) \, d\xi dt 
$$
by Strichartz, using the fact that at least two of the $|\xi_i|$ are $\ge cN$ , 
by
\begin{eqnarray*}
\lefteqn{ c N_{max} N^{-1} \|u_1\|_{L^4_{xt}} \|v_2\|_{L^4_{xt}} 
\|u_3\|_{L^4_{xt}} \|v_4\|_{L^4_{xt}} } \\
& \le & c N_{max} N^{-1} (N_1N_2N_3N_4)^{-1} \|u_1\|_{X^{1,\frac{1}{2}+}} 
\|v_2\|_{Y^{1,\frac{1}{2}+}} 
\|u_3\|_{X^{1,\frac{1}{2}+}} \|v_4\|_{Y^{1,\frac{1}{2}+}} \\ 
& \le & c N^{-2+} N_{max}^{0-} \|u_1\|_{X^{1,\frac{1}{2}+}} 
\|v_2\|_{Y^{1,\frac{1}{2}+}} 
\|u_3\|_{X^{1,\frac{1}{2}+}} \|v_4\|_{Y^{1,\frac{1}{2}+}} \, .
\end{eqnarray*}
\underline{Estimate of $I_{10}$:} We want to show
\begin{eqnarray*}
& & \int_0^{\delta} \int_* 
\left|\frac{m(\xi_4+\xi_5+\xi_6)-m(\xi_4)m(\xi_5)m(\xi_6)}{m(\xi_4)m(\xi_5)m(\xi
_6)}\right| \prod_{i=1}^6 
\widehat{u_i}(\xi_i,t) d\xi dt \\
& & \quad \le c N^{-3+} \prod_{i=1}^6 \|u_i\|_{X^{1,\frac{1}{2}+}} \, . 
\end{eqnarray*}
\underline{Case 1:} At least three of the $|\xi_i|$ are $\ge cN$ . \\
The multiplier is bounded by $c(\frac{N_{max}}{N})^{\frac{3}{2}}$ , so that an 
application of Strichartz 
gives an integral bound
\begin{eqnarray*}
c(\frac{N_{max}}{N})^{\frac{3}{2}} \prod_{i=1}^6 \|u_i\|_{L^6_{xt}} & \le & 
c(\frac{N_{max}}{N})^{\frac{3}{2}} \prod_{i=1}^6 N_i^{-1} \prod_{i=1}^6 
\|u_i\|_{X^{1,\frac{1}{2}+}} \\
& \le & c N^{-3+} N_{max}^{0-} \prod_{i=1}^6 \|u_i\|_{X^{1,\frac{1}{2}+}} \, .
\end{eqnarray*}
\underline{Case 2:} Exactly two of the $|\xi_i|$ are $\ge cN$ , the others $<< 
N$ , e.g. 
$|\xi_5|,|\xi_6|\ge cN$ . \\
Then the multiplier is bounded by 
$c(\frac{N_5}{N})^{\frac{1}{2}}(\frac{N_6}{N})^{\frac{1}{2}}$, and the 
integral, using (\ref{c}), by
\begin{eqnarray*}
& & c(\frac{N_5}{N})^{\frac{1}{2}} (\frac{N_6}{N})^{\frac{1}{2}} \|u_1 
u_5\|_{L^2_{xt}} \|u_2 
u_6\|_{L^2_{xt}} \|u_3\|_{L^{\infty}_{xt}} \|u_4\|_{L^{\infty}_{xt}} \\
& & \quad \le c(\frac{N_5}{N})^{\frac{1}{2}} (\frac{N_6}{N})^{\frac{1}{2}} 
\|u_1\|_{X^{0,\frac{1}{2}+}} 
\|D_x^{-\frac{1}{2}}u_5\|_{X^{0,\frac{1}{2}+}} \|u_2\|_{X^{0,\frac{1}{2}+}} 
\|D_x^{-\frac{1}{2}}u_6\|_{X^{0,\frac{1}{2}+}} \times \\
& & \hspace{20em} \times \|u_3\|_{L^{\infty}_t H^{\frac{1}{2}+}_x} 
\|u_4\|_{L^{\infty}_t 
H^{\frac{1}{2}+}_x} \\
& & \quad \le c(\frac{N_5}{N})^{\frac{1}{2}} (\frac{N_6}{N})^{\frac{1}{2}} 
N_5^{-\frac{3}{2}} 
N_6^{-\frac{3}{2}} \prod_{i=1}^6 \|u_i\|_{X^{1,\frac{1}{2}+}} \\
& & \quad \le c N^{-3+} N_{max}^{0-} \prod_{i=1}^6 \|u_i\|_{X^{1,\frac{1}{2}+}} 
\, . 
\end{eqnarray*}
\underline{Estimate of $I_{11}$ :} We have to show
\begin{eqnarray*}
& & \int_0^{\delta} \int_* 
\left|\frac{m(\xi_1+\xi_2+\xi_3)-m(\xi_1)m(\xi_2)m(\xi_3)}{m(\xi_1)m(\xi_2)m(\xi
_3)}\right| \prod_{i=1}^4 
\widehat{u_i}(\xi_i,t) \widehat{v_5}(\xi_5,t) \, d\xi dt \\
& & \quad \le c N^{-3+} \prod_{i=1}^4 \|u_i\|_{X^{1,\frac{1}{2}+}} 
\|v_5\|_{Y^{1,\frac{1}{2}+}} \, .
\end{eqnarray*}
\underline{Case 1:} At least three of the $|\xi_i|$ are $\ge cN$ . \\
In this case we bound the multiplier by $c(\frac{N_{max}}{N})^{\frac{3}{2}}$ and 
use Strichartz to control 
the integral by
\begin{eqnarray*}
c(\frac{N_{max}}{N})^{\frac{3}{2}} \prod_{i=1}^4 \|u_i\|_{L^5_{xt}} 
\|v_5\|_{L^5_{xt}} & \le & 
c(\frac{N_{max}}{N})^{\frac{3}{2}} \prod_{i=1}^5 N_i^{-1} \prod_{i=1}^4 
\|u_i\|_{X^{1,\frac{1}{2}+}} 
\|v_5\|_{Y^{1,\frac{1}{2}+}} \\
& \le & c N^{-3+} N_{max}^{0-} \prod_{i=1}^4 \|u_i\|_{X^{1,\frac{1}{2}+}} 
\|v_5\|_{Y^{1,\frac{1}{2}+}} \, .
\end{eqnarray*}
\underline{Case 2:} Exactly two of the $|\xi_i|$ are $\ge cN$ , the others $<< 
N$ . \\
a. $|\xi_1|,|\xi_2| \ge cN$ . \\
The multiplier is bounded by 
$c(\frac{N_1}{N})^{\frac{1}{2}}(\frac{N_2}{N})^{\frac{1}{2}}$, and the 
integral, using (\ref{c}), by
\begin{eqnarray*}
& & \hspace{-2em}c(\frac{N_1}{N})^{\frac{1}{2}}(\frac{N_2}{N})^{\frac{1}{2}} 
\|u_1 u_3\|_{L^2_{xt}} \|u_2 
u_4\|_{L^2_{xt}} \|v_5\|_{L^{\infty}_{xt}} \\
& &\hspace{-2em} \le c(\frac{N_1}{N})^{\frac{1}{2}}(\frac{N_2}{N})^{\frac{1}{2}}
 \|D_x^{-\frac{1}{2}} u_1\|_{X^{0,\frac{1}{2}+}} \|u_3\|_{X^{0,\frac{1}{2}+}}
\|D_x^{-\frac{1}{2}} u_2\|_{X^{0,\frac{1}{2}+}} \|u_4\|_{X^{0,\frac{1}{2}+}} 
\|v_5\|_{L^{\infty}_t 
H^{\frac{1}{2}+}_x} \\
& & \hspace{-2em} \le 
c(\frac{N_1}{N})^{\frac{1}{2}}(\frac{N_2}{N})^{\frac{1}{2}} N_1^{-\frac{3}{2}} 
N_2^{-\frac{3}{2}} \prod_{i=1}^4 \|u_i\|_{X^{1,\frac{1}{2}+}} 
\|v_5\|_{Y^{1,\frac{1}{2}+}} \\
& & \hspace{-2em} \le c N^{-3+} N_{max}^{0-} \prod_{i=1}^4 
\|u_i\|_{X^{1,\frac{1}{2}+}} 
\|v_5\|_{Y^{1,\frac{1}{2}+}} \, .
\end{eqnarray*}
b. $|\xi_1|,|\xi_5| \ge cN$ . \\
Estimate the multiplier crudely by a constant and use Strichartz and Lemma 
\ref{Lemma 3} to bound the 
integral by
\begin{eqnarray*}
& & c \prod_{i=1}^3 \|u_i\|_{L^6_{xt}} \|u_4 v_5\|_{L^2_{xt}} \le c 
\prod_{i=1}^4 
\|u_i\|_{X^{0,\frac{1}{2}+}} \|D_x^{-1}v_5\|_{Y^{0,\frac{1}{2}+}} \\
& & \le c N_1^{-1} N_5^{-2} \prod_{i=1}^4 \|u_i\|_{X^{1,\frac{1}{2}+}} 
\|v_5\|_{Y^{1,\frac{1}{2}+}}  \le  c 
N^{-3+} N_{max}^{0-}  \prod_{i=1}^4 \|u_i\|_{X^{1,\frac{1}{2}+}} 
\|v_5\|_{Y^{1,\frac{1}{2}+}} \, .
\end{eqnarray*}
The remaining cases are similar.\\
\underline{Estimate of $I_{12}$ :} similarly as $I_{11}$ .
\section{Estimates for the modified L - functional}
We also need control over the increment of the modified L - functional.
\begin{prop}
\label{Proposition 5.1}
If $(u,v)$ is a solution of (\ref{1.1}),(\ref{1.2}),(\ref{1.3}) on $[0,\delta]$ 
in the sense of Proposition 
\ref{Proposition 1}, then the following estimate holds for $N \ge 1$ and $s > 
1/2$ :
\begin{eqnarray*}
& & |L(Iu(\delta),Iv(\delta))-L(Iu(0),Iv(0))| \\
& & \quad \le c\left[ N^{-2+} \delta^{\frac{1}{2}-} 
\left(\|Iu\|_{X^{1,\frac{1}{2}+}}^3 + 
\|Iv\|_{Y^{1,\frac{1}{2}+}}^3 \right) + N^{-3+} \|Iu\|_{X^{1,\frac{1}{2}+}}^4 
\right] \, .
\end{eqnarray*}
\end{prop}
{\bf Proof:} We use (\ref{3.2'}) and argue similarly as in the previous 
proposition.\\
\underline{Estimate of $J_1$ :} We have to show
\begin{eqnarray*}
& & \int_0^{\delta} 
\int_*\left|\frac{m(\xi_1+\xi_2)-m(\xi_1)m(\xi_2)}{m(\xi_1)m(\xi_2)} \right| 
\widehat{v_1}(\xi_1,t) |\xi_2| \widehat{v_2}(\xi_2,t) \widehat{v_3}(\xi_3,t) \, 
d\xi dt \\
& & \quad \le cN^{-2+} \delta^{\frac{1}{2}-} \prod_{i=1}^3 
\|v_i\|_{Y^{1,\frac{1}{2}+}}  \, .
\end{eqnarray*}
\underline{Case 1:} $|\xi_1| << |\xi_2| \sim |\xi_3| \ge cN $ . \\
a. $|\xi_1| \le N$ . \\
The multiplier is controlled by $ c |\frac{(\nabla m)(\xi_2) \xi_1}{m(\xi_2)}| 
\le c \frac{N_1}{N_2} $ and 
the integral, using (\ref{g}), by
\begin{eqnarray*}
& & c N_1 N_2^{-1} \|D_x^{\frac{1}{2}} I_-^{\frac{1}{2}} (v_1, v_2)\|_{L^2_{xt}} 
\|v_3\|_{L^2_{xt}}  \le 
c N_1 N_2^{-1} \|v_1\|_{Y^{0,\frac{1}{2}+}} \|v_2\|_{Y^{0,\frac{1}{2}+}} 
\|v_3\|_{L^2_{xt}} \\
& & \le c N_1 N_2^{-1} (N_1N_2N_3)^{-1} \delta^{\frac{1}{2}-} \prod_{i=1}^3 
\|v_i\|_{Y^{1,\frac{1}{2}+}} 
\le c N^{-3+} \delta^{\frac{1}{2}-} N_{max}^{0-} \prod_{i=1}^3 
\|v_i\|_{Y^{1,\frac{1}{2}+}} \, .
\end{eqnarray*}
b. $|\xi_1| \ge N$ . \\
Using the multiplier bound $c(\frac{N_1}{N})^{\frac{1}{2}}$ and estimating as in 
case a we get the same.\\
Similarly we treat the case $|\xi_1| >> |\xi_2| $ . \\
\underline{Case 2:} $|\xi_1| \sim |\xi_2| \ge cN $ , $ |\xi_3| << 
|\xi_1|,|\xi_2|$ . \\
This gives the bound, using (\ref{g}):
\begin{eqnarray*}
& & c N_1 N^{-1} \|v_1\|_{L^2_{xt}} 
\|D_x^{\frac{1}{2}}I_-^{\frac{1}{2}}(v_2,v_3)\|_{L^2_{xt}} \le c N_1 
N^{-1} \|v_1\|_{L^2_{xt}} \|v_2\|_{Y^{0,\frac{1}{2}+}} 
\|v_3\|_{Y^{0,\frac{1}{2}+}} \\
& & \le c N_1 N^{-1} N_1^{-1} \delta^{\frac{1}{2}-} N_2^{-1} \prod_{i=1}^3 
\|v_i\|_{Y^{1,\frac{1}{2}+}} \le 
c N^{-2+} N_{max}^{0-} \delta^{\frac{1}{2}-} \prod_{i=1}^3 
\|v_i\|_{Y^{1,\frac{1}{2}+}} \, .
\end{eqnarray*}
\underline{Case 3:} $|\xi_1| \sim |\xi_2| \sim |\xi_3| \ge cN $ . \\
By Strichartz we get the bound
\begin{eqnarray*}
& & c N_1 N^{-1} \|v_1\|_{L^2_{xt}} \|D_x v_2\|_{L^4_{xt}} \|v_3\|_{L^4_{xt}} 
\le c N_1 N^{-1} 
\|v_1\|_{L^2_{xt}} \|v_2\|_{Y^{1,\frac{1}{2}+}} \|v_3\|_{Y^{0,\frac{1}{2}+}} \\
& & \le c N_1 N^{-1} \delta^{\frac{1}{2}-} N_1^{-1} N_3^{-1} \prod_{i=1}^3 
\|v_i\|_{Y^{1,\frac{1}{2}+}} \le 
c N^{-2+} N_{max}^{0-} \delta^{\frac{1}{2}-} \prod_{i=1}^3 
\|v_i\|_{Y^{1,\frac{1}{2}+}} \, .
\end{eqnarray*} 
\underline{Estimate of $J_2$ :} We want to show
\begin{eqnarray*}
& & \int_0^{\delta} \int_* 
\left|\frac{m(\xi_1+\xi_2)-m(\xi_1)m(\xi_2)}{m(\xi_1)m(\xi_2)}\right| 
|\xi_1+\xi_2| \widehat{u_1}(\xi_1,t)  \widehat{u_2}(\xi_2,t)  
\widehat{v_3}(\xi_3,t) \, d\xi dt \\
& & \quad \le c N^{-2+} \delta^{\frac{1}{2}-} \|u_1\|_{X^{1,\frac{1}{2}+}} 
\|u_2\|_{X^{1,\frac{1}{2}+}}
\|v_3\|_{Y^{1,\frac{1}{2}+}} \, .
\end{eqnarray*}
\underline{Case 1:} $|\xi_3| \sim|\xi_1| \ge cN $ , $ |\xi_2| << |\xi_1|,|\xi_3| 
$ . \\
We estimate the multiplier by $c\langle(\frac{N_2}{N})^{\frac{1}{2}}\rangle$ and 
the rest of the integral 
using Lemma \ref{Lemma 3} by
\begin{eqnarray*}
c \|D_x u_1\|_{L^2_{xt}} \|u_2v_3\|_{L^2_{xt}} & \le & c \|D_x u_1\|_{L^2_{xt}} 
\|u_2\|_{X^{0,\frac{1}{2}+}} \|D_x^{-1} v_3\|_{Y^{0,\frac{1}{2}+}} \\
& \le & c N_3^{-2} \delta^{\frac{1}{2}-} \|u_1\|_{X^{1,\frac{1}{2}+}} 
\|u_2\|_{X^{1,\frac{1}{2}+}} 
\|v_3\|_{Y^{1,\frac{1}{2}+}} \, .
\end{eqnarray*}
This gives the desired estimate.\\
\underline{Case 2:} $|\xi_1| \sim |\xi_2| \ge cN$ . \\
By Strichartz we get the bound
\begin{eqnarray*}
\lefteqn{c N_2 N^{-1} \|u_1\|_{L^4_{xt}} \|u_2\|_{L^4_{xt}} \|D_x 
v_3\|_{L^2_{xt}} } \\
 & \le & c N_2 N^{-1} N_1^{-1} N_2^{-1} \delta^{\frac{1}{2}-} 
\|u_1\|_{X^{1,\frac{1}{2}+}}  
\|u_2\|_{X^{1,\frac{1}{2}+}}  \|v_3\|_{Y^{1,\frac{1}{2}+}} \\
& \le & c N^{-2+} N_{max}^{0-} \delta^{\frac{1}{2}-} 
\|u_1\|_{X^{1,\frac{1}{2}+}}  
\|u_2\|_{X^{1,\frac{1}{2}+}}  \|v_3\|_{Y^{1,\frac{1}{2}+}} \, .
\end{eqnarray*}
\underline{Estimate of $J_3$ :} We have to show
\begin{eqnarray*}
& & \int_0^{\delta} \int_* 
\left|\frac{m(\xi_1+\xi_2)-m(\xi_1)m(\xi_2)}{m(\xi_1)m(\xi_2)}\right| 
\widehat{u_1}(\xi_1,t)  \widehat{v_2}(\xi_2,t) |\xi_3| \widehat{u_3}(\xi_3,t) \, 
d\xi dt \\
& & \quad \le c N^{-2+} \delta^{\frac{1}{2}-} \|u_1\|_{X^{1,\frac{1}{2}+}} 
\|v_2\|_{Y^{1,\frac{1}{2}+}}
\|u_3\|_{X^{1,\frac{1}{2}+}} \, .
\end{eqnarray*}
A typical case is $|\xi_3| \sim |\xi_1| \ge cN $ , $ |\xi_2| << |\xi_1|,|\xi_3| 
$ . Estimating the 
multiplier by $c\frac{N_2}{N_1}$ , if $|\xi_2|\le N$ , and by 
$c(\frac{N_2}{N})^{\frac{1}{2}}$ , if 
$|\xi_2| \ge N$ , we get by Strichartz the folllowing bound for the rest of the 
integral
$$\|u_1\|_{L^4_{xt}} \|v_2\|_{L^2_{xt}} \|D_x u_3 \|_{L^4_{xt}} \le N_1^{-1} 
N_2^{-1} \delta^{\frac{1}{2}-}
\|u_1\|_{X^{1,\frac{1}{2}+}} \|v_2\|_{Y^{1,\frac{1}{2}+}} \|u_3 
\|_{X^{1,\frac{1}{2}+}} \, , 
$$
which gives the claimed estimate. The other cases are treated similarly.\\
\underline{Estimate of $J_4$ :} The desired estimate is
\begin{eqnarray*}
& & \int_0^{\delta} \int_* 
\left|\frac{m(\xi_1+\xi_2+\xi_3)-m(\xi_1)m(\xi_2)m(\xi_3)}{m(\xi_1)m(\xi_2)m(\xi
_3)}\right| \prod_{i=1}^4 
\widehat{u_i}(\xi_i,t)|\xi_4| \, d\xi dt \\ 
& & \quad \le c N^{-3+} \prod_{i=1}^4 \|u_i\|_{X^{1,\frac{1}{2}+}} \, .
\end{eqnarray*}
\underline{Case 1:} $N_1,N_2,N_3 \ge cN$ . \\
Estimate the multiplier by $\, 
c(\frac{N_1}{N}\frac{N_2}{N}\frac{N_3}{N})^{\frac{1}{2}} \, $ and the rest 
of the integral using \\ 
Strichartz by
$$ \prod_{i=1}^3 \|u_i\|_{L^4_{xt}} \|D_x u_4\|_{L^4_{xt}} \le c 
(N_1N_2N_3)^{-1} \prod_{i=1}^4 
\|u_i\|_{X^{1,\frac{1}{2}+}} \, . $$
This gives the desired bound. \\
\underline{Case 2:} $N_1 \sim N_2 \ge cN $ , $ N_3,N_4 << N_1,N_2 $ . \\
This gives the bound by (\ref{c}) :
\begin{eqnarray*}
& & \hspace{-2em} c (\frac{N_1}{N})^{\frac{1}{2}} (\frac{N_2}{N})^{\frac{1}{2}} 
\langle(\frac{N_3}{N})^{\frac{1}{2}}\rangle \|u_1 u_3\|_{L^2_{xt}} \|u_2 D_x 
u_4\|_{L^2_{xt}} \\
& & \hspace{-2em} \le c (\frac{N_1}{N})^{\frac{1}{2}} 
(\frac{N_2}{N})^{\frac{1}{2}} 
\langle(\frac{N_3}{N})^{\frac{1}{2}}\rangle \|D_x^{-\frac{1}{2}} 
u_1\|_{X^{0,\frac{1}{2}+}} 
\|u_3\|_{X^{0,\frac{1}{2}+}} \|D_x^{-\frac{1}{2}} u_2\|_{X^{0,\frac{1}{2}+}} 
\|D_x 
u_4\|_{X^{0,\frac{1}{2}+}} \\
& &  \hspace{-2em}\le c (\frac{N_1}{N})^{\frac{1}{2}} 
(\frac{N_2}{N})^{\frac{1}{2}} 
\langle(\frac{N_3}{N})^{\frac{1}{2}}\rangle (N_1N_2)^{-\frac{3}{2}} \langle 
N_3\rangle^{-1} \prod_{i=1}^4 
\|u_i\|_{X^{1,\frac{1}{2}+}} \\
& & \hspace{-2em}\le c N^{-3+} N_{max}^{0-} \prod_{i=1}^4 
\|u_i\|_{X^{1,\frac{1}{2}+}} \, .
\end{eqnarray*}
\underline{Case 3:} $N_1 \sim N_4 \ge cN $ , $ N_2,N_3 << N_1,N_4 $ . \\
a. $N_2,N_3 \le N$ . \\
The multiplier is bounded by $ c|\frac{(\nabla m)(\xi_1)}{m(\xi_1)} 
(\xi_2+\xi_3)| \le c 
\frac{N_2+N_3}{N_1} $ , and thus by (\ref{c}) we get the bound
\begin{eqnarray*}
\lefteqn{ c \frac{N_2+N_3}{N_1} \|u_1 u_2\|_{L^2_{xt}} \|u_3 D_x 
u_4\|_{L^2_{xt}} } \\
& \le & c \frac{N_2+N_3}{N_1} \|u_1\|_{X^{-\frac{1}{2},\frac{1}{2}+}} 
\|u_2\|_{X^{0,\frac{1}{2}+}} 
\|u_3\|_{X^{0,\frac{1}{2}+}} \|D_x^{\frac{1}{2}} u_4\|_{X^{0,\frac{1}{2}+}} \\
& \le & \frac{N_2+N_3}{N_1} N_1^{-\frac{3}{2}} \langle N_2 \rangle^{-1} \langle 
N_3 \rangle^{-1} 
N_4^{-\frac{1}{2}} \prod_{i=1}^4 \|u_i\|_{X^{1,\frac{1}{2}+}} \\
& \le & c N^{-3+} N_{max}^{0-} \prod_{i=1}^4 \|u_i\|_{X^{1,\frac{1}{2}+}} \, .
\end{eqnarray*}
b. $N_2 \ge N$ (similarly $N_3 \ge N$) . \\
The estimate is similar to case a, but the multiplier bound is 
$c(\frac{N_1}{N})^{\frac{1}{2}}(\frac{N_2}{N})^{\frac{1}{2}}\langle(\frac{N_3}{N
})^{\frac{1}{2}}\rangle.$ 
This leads to 
the same bound as in a. \\
\underline{Case 4:} $N_1 \sim N_2 \sim N_4 \ge cN $ , $ N_3 << N_1,N_2,N_4 $ . 
\\
Because of $\sum_{i=1}^4 \xi_i = 0$, two of the large frequencies have different 
sign, say, $\xi_2$ and 
$\xi_4$. Thus $|\xi_4|^{\frac{1}{2}} \le |\xi_2-\xi_4|^{\frac{1}{2}}$, and we 
get the bound for the 
integral, using (\ref{h}) and (\ref{c}):
\begin{eqnarray*}
\lefteqn{ c (\frac{N_1}{N})^{\frac{1}{2}}(\frac{N_2}{N})^{\frac{1}{2}}\langle 
(\frac{N_3}{N})^{\frac{1}{2}} 
\rangle \|u_2 D_x u_4\|_{L^2_{xt}} \|u_1 u_3\|_{L^2_{xt}} } \\
& \le & c (\frac{N_1}{N})^{\frac{1}{2}}(\frac{N_2}{N})^{\frac{1}{2}}\langle 
(\frac{N_3}{N})^{\frac{1}{2}} 
\rangle \|I_-^{\frac{1}{2}}(u_2, D_x^{\frac{1}{2}} u_4)\|_{L^2_{xt}} \|u_1 
u_3\|_{L^2_{xt}} \\
& \le & c 
(\frac{N_1}{N})^{\frac{1}{2}}(\frac{N_2}{N})^{\frac{1}{2}}\langle(\frac{N_3}{N})
^{\frac{1}{2}}\rangle 
\|u_2\|_{X^{0,\frac{1}{2}+}} \|D_x^{\frac{1}{2}}u_4\|_{X^{0,\frac{1}{2}+}} 
\|D_x^{-\frac{1}{2}}u_1\|_{X^{0,\frac{1}{2}+}} \|u_3\|_{X^{0,\frac{1}{2}+}} \\
& \le & c 
(\frac{N_1}{N})^{\frac{1}{2}}(\frac{N_2}{N})^{\frac{1}{2}}\langle(\frac{N_3}{N})
^{\frac{1}{2}}\rangle 
N_2^{-1} N_4^{-\frac{1}{2}} N_1^{-\frac{3}{2}} \langle N_3 \rangle^{-1} 
\prod_{i=1}^4 
\|u_i\|_{X^{1,\frac{1}{2}+}} \\
& \le & c N^{-3+} N_{max}^{0-} \prod_{i=1}^4 \|u_i\|_{X^{1,\frac{1}{2}+}} \, .
\end{eqnarray*}
\section{The global existence result}
\begin{theorem}
Let $ 1 > s > 3/5 $ , if $\beta =0$ , and $ 1 > s > 2/3 $ , if $ \beta \neq 0 $ 
, and $\alpha \gamma > 0 $ 
. The system (\ref{1.1}),(\ref{1.2}),(\ref{1.3}) has a unique global solution 
for data $(u_0,v_0) \in 
H^s({\bf R}) \times H^s({\bf R}) $ . More precisely, for any $T>0$ there exists 
$b> 1/2$ and a unique 
solution $(u,v) \in X^{s,b}[0,T] \times Y^{s,b}[0,T]$ with $(u,v) \in 
C^0([0,T],H^s({\bf R})\times H^s({\bf 
R})) $ . 
\end{theorem}
{\bf Proof:} The data satisfy the estimates
$ \|Iu_0\|_{H^1}^2 + \|Iv_0\|_{H^1}^2 \le c N^{2(1-s)}$ and $ \|Iu_0\|_{L^2}^2 + 
\|Iv_0\|_{L^2}^2 \le c $ , 
especially $ \|Iv_0\|_{L^2}^2 \le c N^{1-s} $ . These bounds imply by 
(\ref{3.1}) and (\ref{3.4}): 
$|L(Iu_0,Iv_0)| \le \overline{c} N^{1-s} $ and $|E(Iu_0,Iv_0)| \le \overline{c} 
N^{2(1-s)} $ , and any such 
bounds for $L$ and $E$ imply by (\ref{3.5}) and (\ref{3.6}): $\|Iu_0\|_{L^2}^2 
\le M^2 $ , $ 
\|Iv_0\|_{L^2}^2 \le \widehat{c} N^{1-s} $ , $ \|Iu_0\|_{H^1}^2 + 
\|Iv_0\|_{H^1}^2 \le \widehat{c} 
N^{2(1-s)} $ with $\widehat{c}=\widehat{c}(\overline{c})$ . \\
We use our local existence theorem on $[0,\delta]$, where $\delta \sim 
N^{-4(1-s)+}$ , if $\beta \neq 0 $ , 
and $ \delta \sim N^{-3(1-s)+} $ , if $\beta = 0$ , and conclude
\begin{equation}
\label{6.1}
\|Iu\|_{X^{1,\frac{1}{2}+}[0,\delta]} + \|Iv\|_{Y^{1,\frac{1}{2}+}[0,\delta]} 
\le c_0 (\|Iu_0\|_{H^1} + 
\|Iv_0\|_{H^1}) \le c_1 N^{1-s}\, ,
\end{equation}
where $c_1 = c_1(\overline{c},M)$ . In order to reapply the local existence 
theorem with intervals of equal 
length we need a uniform bound of the $H^1$ - norms of the solution at time 
$t=\delta$ and $t=2\delta$ etc. 
This follows from uniform control over $|E|$ and $|L|$ by (\ref{3.6}). The 
increment of $E$ is controlled 
by Proposition \ref{Proposition 4.1} and (\ref{6.1}) as follows:
\begin{eqnarray*}
\lefteqn{ |E(Iu(\delta),Iv(\delta)) - E(Iu_0,Iv_0)| } \\
& \le & c_2[(N^{-1+} \delta^{\frac{1}{2}-} + N^{-\frac{7}{4}+})N^{3(1-s)} + 
N^{-2+} N^{4(1-s)} + N^{-3+} 
N^{6(1-s)}] \, , 
\end{eqnarray*}
where $ c_2 = c_2(\overline{c},M)$ . \\
The number of iteration steps to reach the given time $T$ is $T \delta^{-1}$ . 
This means that in order to 
give a uniform bound of the energy of the iterated solutions by $2\overline{c} 
N^{2(1-s)}$ , the following 
condition has to be fulfilled:
\begin{equation}
\label{6.2}
c_2[(N^{-1+} \delta^{\frac{1}{2}-} + N^{-\frac{7}{4}+})N^{3(1-s)} + N^{-2+} 
N^{4(1-s)} + N^{-3+}N^{6(1-s)}] 
T \delta^{-1} < \overline{c} N^{2(1-s)} \, ,
\end{equation}
where $c_2 = c_2(2\overline{c},2M)$ (recall that the initial energy is bounded 
by $\overline{c} 
N^{2(1-s)}$).
Similarly, the increment of $L$ is controlled by Proposition \ref{Proposition 
5.1} and (\ref{6.1}):
$$|L(Iu(\delta),Iv(\delta)) - L(Iu_0,Iv_0)| \le c_2[N^{-2+} 
\delta^{\frac{1}{2}-} N^{3(1-s)} + N^{-3+} 
N^{4(1-s)}] \, . $$
Thus, similarly as for $E$, in order to give a uniform bound of $L$ by 
$2\overline{c}N^{1-s}$ , the 
following condition has to be fulfilled:
\begin{equation}
\label{6.3}
c_2(N^{-2+} \delta^{\frac{1}{2}-} N^{3(1-s)} + N^{-3+} N^{4(1-s)}) T \delta^{-1} 
< \overline{c} N^{1-s} \, 
. 
\end{equation}
If the inequalities (\ref{6.2}) and (\ref{6.3}) are satisfied, the uniform 
control of $|E|$ and $|L|$ 
implies by (\ref{3.5}) and(\ref{3.6}) uniform control 
$$ \|v(t)\|_{L^2} \le c N^{\frac{1-s}{2}} \quad \mbox{and} \quad \|u(t)\|_{H^1} 
+ \|v(t)\|_{H^1} \le c  
N^{1-s} \, . $$
Now, using the definiton of $\delta$ above, (\ref{6.2}) can be fulfilled for a 
sufficiently large $N$, 
provided the following conditons hold: \\
a. in the case $\beta \neq 0$ : $-1-2(1-s)+3(1-s)+4(1-s)<2(1-s) \Leftrightarrow 
s > 2/3$, and 
$-\frac{7}{4}+3(1-s)+4(1-s)<2(1-s) \Leftrightarrow s > 13/20$, and 
$-2+4(1-s)+4(1-s)<2(1-s) \Leftrightarrow 
s > 2/3$, and $-3+6(1-s)+4(1-s)<2(1-s) \Leftrightarrow s > 5/8$.\\
b. in the case $\beta=0$ : $-1-\frac{3}{2}(1-s)+3(1-s)+3(1-s)<2(1-s) 
\Leftrightarrow s> 3/5$ and 
$-\frac{7}{4}+3(1-s)+3(1-s)<2(1-s) \Leftrightarrow s > 9/16$ , and 
$-2+4(1-s)+3(1-s)<2(1-s) \Leftrightarrow 
s> 3/5$ , and $-3+6(1-s)+3(1-s)<2(1-s) \Leftrightarrow s > 4/7$ . \\
Similarly, (\ref{6.3}) is fulfilled for $N$ sufficiently large, provided \\
a. in the case $\beta \neq 0$ : $-2-2(1-s)+3(1-s)+4(1-s)< 1-s \Leftrightarrow s 
> 1/2$ , and 
$-3+4(1-s)+4(1-s)<1-s \Leftrightarrow s > 4/7$ . \\
b. in the case $\beta =0$ : $-2-\frac{3}{2}(1-s)+3(1-s)+3(1-s)<1-s 
\Leftrightarrow s > 3/7$ , and 
$-3+4(1-s)+3(1-s)<1-s \Leftrightarrow s > 1/2$ .
 
\end{document}